\documentclass[11pt]{article}
\usepackage{color}
\usepackage{amsmath}
\usepackage{amssymb}
\usepackage{latexsym}
\usepackage{comment}

\newtheorem{assumption}{Assumption}

\def\qed{ \ \vrule width.2cm height.2cm depth0cm\smallskip}
\newenvironment{proof}{\noindent {\bf Proof.\/}}{$\qed$\vskip 0.1in}

\newcommand{\la}{\langle}
\newcommand{\ra}{\rangle}

\newcommand{\ba}{\begin{array}}
\newcommand{\ea}{\end{array}}
\newcommand{\be}{\begin{equation}}
\newcommand{\ee}{\end{equation}}
\newcommand{\bea}{\begin{eqnarray}}
\newcommand{\eea}{\end{eqnarray}}
\newcommand{\beaa}{\begin{eqnarray*}}
\newcommand{\eeaa}{\end{eqnarray*}}

\def\dbC{\mathbb{C}}
\def\dbD{\mathbb{D}}
\def\dbE{\mathbb{E}}
\def\dbF{\mathbb{F}}

\def\dbK{\mathbb{K}}
\def\dbL{\mathbb{L}}

\def\dbP{\mathbb{P}}
\def\dbR{\mathbb{R}}
\def\dbS{\mathbb{S}}
\def\dbQ{\mathbb{Q}}

\def\a{\alpha}
\def\b{\beta}
\def\g{\gamma}
\def\d{\delta}
\def\e{\varepsilon}

\def\l{\lambda}

\def\si{\sigma}
\def\t{\tau}
\def\f{\varphi}
\def\th{\theta}
\def\o{\omega}

\def\L{\Lambda}

\def\O{\Omega}
\def\cA{{\cal A}}

\def\cC{{\cal C}}
\def\cD{{\cal D}}
\def\cE{{\cal E}}
\def\cF{{\cal F}}

\def\cH{{\cal H}}

\def\cK{{\cal K}}
\def\cL{{\cal L}}

\def\cP{{\cal P}}

\def\cS{{\cal S}}
\def\cT{{\cal T}}
\def\cU{{\cal U}}

\def\cX{{\cal X}}
\def\cY{{\cal Y}}
\def\cZ{{\cal Z}}

\def\ch{\textsc{h}}

\def\no{\noindent}

\def\ms{\medskip}

\def\q{\quad}

\def\pa{\partial}
\def\cd{\cdot}
\def\cds{\cdots}

\def\tr{\hbox{\rm tr}}

\def\qed{ \hfill \vrule width.25cm height.25cm depth0cm\smallskip}

\newcommand{\basa}{\begin{assumption}}
\newcommand{\easa}{\end{assumption}}

\newcommand{\bas}{\begin{assum}}
\newcommand{\eas}{\end{assum}}

\def\pa{\partial}

 \def\cd{\cdot}
\def\cds{\cdots}

\def\sgn{\hbox{\rm sgn$\,$}}

\def\tr{\hbox{\rm tr$\,$}}

\def\dis{\displaystyle}

\def\wh{\widehat}

\def\cad{c\`{a}dl\`{a}g}

\def\1{{\bf 1}}

\def\:{\!:\!}
\def\reff#1{{\rm(\ref{#1})}}
\def \proof{{\noindent \bf Proof\quad}}

\def \dbf{{\mathbf{d}}}

\def \Usup{\overline{\cU}}
\def \Usub{\underline{\cU}}

at 9pt

\begin{document}

\newtheorem{thm}{Theorem}[section]
\newtheorem{lem}[thm]{Lemma}
\newtheorem{cor}[thm]{Corollary}
\newtheorem{prop}[thm]{Proposition}
\newtheorem{rem}[thm]{Remark}
\newtheorem{eg}[thm]{Example}
\newtheorem{defn}[thm]{Definition}
\newtheorem{assum}[thm]{Assumption}

\renewcommand {\theequation}{\arabic{section}.\arabic{equation}}
\def\thesection{\arabic{section}}

\title{\bf  Viscosity Solutions of Fully Nonlinear Parabolic Path Dependent PDEs: Part I}

\author{Ibrahim  {\sc Ekren}\footnote{University of Southern California, Department of Mathematics, ekren@usc.edu.}   \and 
Nizar {\sc Touzi}\footnote{CMAP, Ecole Polytechnique Paris, nizar.touzi@polytechnique.edu. Research supported by the Chair {\it Financial Risks} of the {\it Risk Foundation} sponsored by Soci\'et\'e
              G\'en\'erale, and
              the Chair {\it Finance and Sustainable Development} sponsored by EDF and Calyon. }
       \and Jianfeng {\sc Zhang}\footnote{University of Southern California, Department of Mathematics, jianfenz@usc.edu. Research supported in part by NSF grant DMS 10-08873.}
}\maketitle

\begin{abstract}
The main objective of this paper and the accompanying one \cite{ETZ2} is to provide a notion of viscosity solutions for fully nonlinear parabolic path-dependent PDEs. Our definition extends our previous work \cite{EKTZ}, focused on the semilinear case, and is crucially based on the nonlinear optimal stopping problem analyzed in \cite{ETZ0}. We prove that our notion of viscosity solutions is consistent with the corresponding notion of classical solutions, and satisfies a stability property and a partial comparison result. The latter is a key step for the wellposedness results established in \cite{ETZ2}. We also show that the value processes of path-dependent stochastic control problems are viscosity solutions of the corresponding path-dependent dynamic programming equations.
\end{abstract}



\vfill\eject

\section{Introduction}
\label{sect-Introduction}
\setcounter{equation}{0}

The objective of this paper is to introduce a notion of viscosity solution of the following fully nonlinear path-dependent partial differential equation:
\bea\label{PPDE-intro}
 -\pa_t u(t,\o)
 - G\big(t,\o, u(t,\o), 
         \pa_\o u(t,\o), \pa^2_{\o\o} u(t,\o)\big)
 =0,
 ~~0\le t< T,~\o\in \O,
\eea
where the unknown $u$ is a progressively measurable process on the canonical space $\O := \{\o\in C([0,T],\dbR^d): \o_0={\bf 0}\}$, and
the nonlinearity $G :  [0,T]\times\O \times  \dbR \times \dbR^d \times \dbS^d \rightarrow \dbR$ is progressively measurable, satisfies convenient Lipschitz and continuity assumptions, and is degenerate elliptic.

The above equation attracted our attention after the point raised by Peng in \cite{Peng-ICM} that this would be an alternative approach to the theory of backward stochastic differential equations, introduced by the seminal paper of Pardoux and Peng \cite{PP1}.

The semilinear case, corresponding to the case where $G$ is linear in the $\partial^2_{\o\o}u-$variable, was addressed in \cite{EKTZ}, where existence and uniqueness results are established for a new notion of viscosity solution. The main difficulty is related to the fact that the canonical space fails to be locally compact, so that many tools from the standard theory of viscosity solutions do not apply to the present context. The main contribution of \cite{EKTZ} is to replace the pointwise extremality in the standard definition of viscosity solutions by the corresponding extremality in the context of an optimal stopping problem under a nonlinear expectation $\cE$. More precisely, we introduce a set of smooth test processes $\varphi$ which are tangent from above or from below to the processes of interest $u$ in the sense of the following nonlinear optimal stopping problems
 \beaa
 \sup_{\tau}\overline{\cE}[(\varphi-u)_{\tau}],
 ~~
 \inf_{\tau}\underline{\cE}[(\varphi-u)_{\tau}],
 &\mbox{where}&
 \overline{\cE}:=\sup_{\dbP\in\cP}\dbE^\dbP,
 ~~\underline{\cE}:=\inf_{\dbP\in\cP}\dbE^\dbP,
 \eeaa
$\tau$ ranges over a convenient set of stopping times, and $\cP$ is a weakly compact collection of probability measures, motivated by a convenient linearization of the nonlinearity $F$. Consequently, in the particular semilinear case of \cite{EKTZ}, the family $\cP$ consists of equivalent probability measures.

In this paper, together with the accompanying ones \cite{ETZ0, ETZ2}, we extend the notion of viscosity solutions to the fully nonlinear case.  As in  \cite{EKTZ}, in order to avoid the local compactness issue of canonical space, we shall still use the optimal stopping problem to define viscosity solutions. However,  in this fully nonlinear context, the family $\cP$ of probability measures consists of nondominated mutually singular measures,  so as to cover all the measures induced by certain linearization of the generator $G$.  The analysis under the corresponding nonlinear expectation $\overline\cE^\cP$ is the major difficulty,  mainly due to the failure of the dominated convergence theorem under $\overline\cE^\dbP$. To overcome this difficulty, one needs some strong regularity of the involved processes which requires rather sophisticated estimates. The corresponding optimal stopping problem is solved in \cite{ETZ0}, and the major result, the comparison principle of viscosity solutions, will be proved in \cite{ETZ2}. 

In this paper we focus on the definition of viscosity solutions and its basic properties. We first prove that our definition of viscosity solutions is consistent with the corresponding notion of classical solutions. Next we show that our viscosity solution satisfies a stability property similar to the finite-dimensional context. Finally, we establish  the partial comparison result, namely for any pair of viscosity subsolution $u^1$ and supersolution $u^2$ with $u^1_T\le u^2_T$ on $\O$, we have $u^1\le u^2$ on $[0, T]\times \O$ whenever either one of them is smooth. This result is crucial for the well-posedness results established in our accompanying paper \cite{ETZ2}.  We remark that Peng \cite{Peng-viscosity} also investigated the comparison principle for fully nonlinear PPDEs by using a  different approach.

We will investigate the connection between our viscosity solution and some other equations in the literature, which will be very helpful for the applications of our results. In particular, we show that the value function of path-dependent stochastic control problems as well as second order backward stochastic differential equations \cite{CSTV,STZ-2BSDE} are naturally viscosity solutions of the corresponding path-dependent partial differential equation. This extends the context of backward stochastic differential equations of \cite{EKTZ}. See also the closely related works \cite{Peng-G, HJPS} in terms of the $G$-expectation. Our PPDE can go beyond stochastic control, see Pham-Zhang \cite{PZ} for an application in stochastic differential games. Moreover, backward stochastic partial differential equations, which can be viewed as the value function of stochastic control with random coefficients (in contrast with path dependent coefficients) can also be viewed as a PPDE. See also \cite{MYZ} and \cite{OSZ} for applications of BSPDEs. 

While the wellposedness of semilinear PPDEs has been achieved in \cite{EKTZ}, the approach there for the comparison principle does not seem to work in fully nonlinear case. We shall revisit the semilinear case by providing a new approach which, modulus all the technicality, will be extended to the fully nonlinear case in \cite{ETZ2}.
Moreover, our context covers first order path-dependent PDEs, which has been studied by Lukoyanov \cite{Lukoyanov} by using compactness arguments.

The rest of the paper is organized as follows. In Section \ref{sect-preliminary}, we introduce the general framework, and define a notion of classical differentiability which is weaker than that of \cite{Dupire}. In Section \ref{sect-PPDE}, we introduce our notion of viscosity solution of fully nonlinear PPDE, and provide various remarks which highlight the analogy with the properties of viscosity solutions in finite dimensional spaces. We prove the consistency with the notion of classical solution. In Section \ref{sect-Example}, we provide several examples and show that natural ones as the value function of path dependent stochastic control problems, or solutions of second order backward stochastic differential equations, are viscosity solutions of the corresponding path-dependent PDEs.
Section \ref{sect-stability} contains our stability and partial comparison results. Section \ref{sect-BSPDE} shows that our framework includes backward stochastic PDE by a convenient augmentation of the canonical space.  Section \ref{sect-semilinear} revisits the semilinear case and provides an alternative and simpler well-posedness argument to that of our previous paper \cite{EKTZ} which will be extended to the fully nonlinear case in our accompanying paper \cite{ETZ2}. Finally, in Section \ref{sect-First} we investigate the first order PPDEs . 

\section{Preliminaries}
\label{sect-preliminary}
\setcounter{equation}{0}

\subsection{The canonical spaces}

Let $\O:= \big\{\o\in C([0,T], \dbR^d): \o_0={\bf 0}\big\}$, the set of continuous paths starting from the origin, $B$ the canonical process, $\dbF$ the natural filtration generated by $B$,  $\dbP_0$ the Wiener measure,  and $\L := [0,T]\times \O$. Here and in the sequel, for notational simplicity,  we use ${\bf 0}$ to denote vectors, matrices, or paths with appropriate dimensions whose components are all equal to $0$. Let $\dbS^d$ denote the set of $d\times d$ symmetric matrices, and
 \beaa
 &x \cd x' := \sum_{i=1}^d x_i x'_i 
 ~~\mbox{for any}~~x, x' \in \dbR^d,
 ~~\g : \g' := \tr[\g\g']
 ~~\mbox{for any}~~\g, \g'\in \dbS^d.
 &
 \eeaa
  We define a seminorm on $\O$ and a pseudometric on $\L$ as follows: for any $(t, \o), ( t', \o') \in\L$, 
\bea\label{rho}
 \|\o\|_{t} 
 := 
 \sup_{0\le s\le t} |\o_s|,
 \q  
 \dbf_\infty\big((t, \o),( t', \o')\big) 
 := 
 |t-t'| + \big\|\o_{.\wedge t} - \o'_{.\wedge t'}\big\|_T.
 \eea
Then $(\O, \|\cd\|_{T})$ is a Banach space and $(\L, \dbf_\infty)$ is a complete pseudometric space. In fact, the subspace  $\{(t, \o_{\cd\wedge t}): (t,\o)\in \L\}$ is a complete metric space under $\dbf_\infty$.  
We shall denote by $\dbL^0(\cF_T)$ and $\dbL^0(\L)$ the collection of all $\cF_T$-measurable random variables and $\dbF$-progressively measurable processes, respectively.  Let $C^0(\L)$ (resp. $UC(\L)$) be the subset of $\dbL^0(\L)$  whose elements are continuous (resp. uniformly continuous) in $(t,\o)$ under $\dbf_\infty$, and $C^0_b(\L)$ (resp. $UC_b(\L)$) be  the subset of $C^0 (\L)$ (resp $UC^0 (\L)$) whose elements are bounded. Finally, $\dbL^0(\L, \dbR^d)$ denote the space of $\dbR^d$-valued processes with entries in $\dbL^0(\L)$, and we define similar notations for the spaces $C^0$, $C^0_b$, $UC$, and $UC_b$.

We next introduce the shifted spaces. Let  $0\le t\le s\le T$.

-  Let $\O^t:= \big\{\o\in C([t,T], \dbR^d): \o_t ={\bf 0}\big\}$ be the shifted canonical space; $B^{t}$ the shifted canonical process on
$\O^t$;   $\dbF^{t}$ the shifted filtration generated by $B^{t}$, $\dbP^t_0$ the Wiener measure on $\O^t$, and $\L^t :=  [t,T]\times \O^t$.

- Define $\|\cd\|^t_s$ on $\O^t$ and $\dbf^t_\infty$ on $\L^t$ in the spirit of (\ref{rho}),  and the sets $\dbL^0(\L^t)$ etc. in an obvious way.   

- For  $\o\in \O^t$ and $\o'\in \O^s$, define the concatenation path $\o\otimes_{s} \o'\in \O^t$ by:
\beaa
(\o\otimes_s \o') (r) := \o_r\1_{[t,s)}(r) + (\o_{s} + \o'_r)\1_{[s, T]}(r),
&\mbox{for all}&
r\in [t,T].
\eeaa

- Let $\xi \in \dbL^0(\cF^t_T)$, and $X\in \dbL^0(\L^t)$. For $(s, \o) \in \L^t$, define $\xi^{s,\o} \in \dbL^0(\cF^s_T)$ and $X^{s,\o}\in \dbL^0(\L^s)$ by:
\beaa
\xi^{s, \o}(\o') :=\xi(\o\otimes_s \o'), \q X^{s, \o}(\o') := X(\o\otimes_s \o'),
&\mbox{for all}&
\o'\in\O^s.
\eeaa

It is clear that, for any $(t,\o) \in \L$ and any $u\in C^0(\L)$, we have $u^{t,\o} \in C^0(\L^t)$. The other spaces introduced before enjoy the same property.

We shall use the following type of regularity, which is slightly stronger than the right continuity of a process $u$  in standard sense (that is, for any fixed $\o$, the mapping $t \mapsto u(t,\o)$ is right continuous). 

\begin{defn}
\label{defn-rightcont}
We say  a process $u \in \dbL^0 (\L)$ is right continuous in $(t,\o)$ under $\dbf_\infty$ if: for any $(t,\o) \in \L$ and any $\e>0$, there exists $\d>0$ such that, for any $(t', \o')\in \L^t$ satisfying $\dbf^t_\infty((t',  \o'), (t, {\bf 0})) \le \d$, we have $|u^{t,\o}(t',\o') - u(t,\o)|\le \e$.
\end{defn} 

\begin{defn}
\label{defn-spaceC0} By $\Usub$, we denote the collection of all 
processes $u \in \dbL^0(\L)$ such that 

- 
$u$ is bounded from above and  right continuous in $(t,\o)$ under $\dbf_\infty$;

- there exists a modulus of continuity function $\rho$ such that for any $(t,\o), (t',\o')\in\L$:
 \bea\label{USC}
 u(t,\o) - u(t',\o') 
 \le 
 \rho\Big(\dbf_\infty\big((t,\o), (t',\o')\big)\Big)
 ~\mbox{whenever}~t\le t'.
 \eea 
By $\Usup$ we denote the set of all processes $u$ such that $-u \in\Usub$.
\end{defn}

\begin{rem}
\label{rem-C0b}
{\rm 
(i) The progressive measurability of $u$  implies that $u(t,\o) = u(t,\o_{\cd\wedge t})$, and it is clear that $\Usub \cap \Usup = UC_b(\L)$. We also recall from \cite{ETZ0} Remark 3.2 
 that Condition \reff{USC} implies that $u$ has left-limits and $u_{t-} \le u_t$ for all $t\in(0,T]$. Moreover, under \reff{USC},  $u$ is right continuous in $(t,\o)$ under $\dbf_\infty$ if and only if it is right continuous in $t$ for every $\o$.
\\
(ii) In finite dimensional case, a continuous function is at least locally uniformly continuous. This is not true anymore in the infinite dimensional case, so it is important to distinguish $C^0(\L)$ and $UC(\L)$ in this paper.
\qed}
\end{rem}

Finally, we denote by $\cT$ the set of $\dbF$-stopping times, and $\cH\subset\cT$ the subset of those hitting times $\ch$ of the form
\bea
\label{cT}
\ch := \inf\{t\ge 0: B_t \in O^c\} \wedge t_0,
\eea
for some $0< t_0\le T$, and some open and convex set $O  \subset \dbR^d$ containing ${\bf 0}$ with $O^c := \dbR^d\setminus O$.
The set $\cH$ will be important for our optimal stopping problem, which is crucial for the comparison and the stability results, see Remark \ref{rem-h}. We note that $\ch = t_0$ when $O=\dbR^d$, and for any $\ch \in \cH$,
\bea
\label{che}
0<\ch_\e \le \ch~~\mbox{for $\e$ small enough, where}~ \ch_\e := \inf\{t\ge 0: |B_t|=\e\} \wedge \e.
\eea
 Moreover,
\beaa
\mbox{$\ch:\O\rightarrow[0,T]$ is lower semicontinuous, and $\ch_1 \wedge \ch_2 \in \cH$ for any $\ch_1, \ch_2\in\cH$}.
\eeaa
Define $\cT^t$ and $\cH^t$ in the same spirit. For any $\t\in \cT$ (resp. $\ch\in\cH$) and any $(t,\o)\in\L$ such that $t<\t(\o)$ (resp. $t<\ch(\o)$), it is clear that $\t^{t,\o}\in \cT^t$ (resp. $\ch^{t,\o}\in\cH^t$).

\subsection{Capacity and nonlinear expectation}

 For every constant $L>0$, we denote by $\cP_L$ the collection of all continuous semimartingale measures $\dbP$ on $\O$ whose drift and diffusion characteristics are bounded by $L$ and $\sqrt{2L}$, respectively.  To be precise, let $\tilde\O:=\O^3$ be an enlarged canonical space, $\tilde B:= (B, A, M)$ be the canonical processes, and $\tilde\o = (\o, a, m)\in \tilde\O$ be the paths. $\dbP\in \cP_L$ means that there exists an extension $\dbQ$ of $\dbP$ on $\tilde\O$ such that: 
\bea
\label{cPL}
\left.\ba{c}
B=A+M,\q \mbox{ $A$ is absolutely continuous, $M$ is a martingale},\\
|\alpha^\dbP|\le L,~~{1\over 2}\tr((\beta^\dbP)^2)\le L,\q \mbox{where}~ \a^\dbP_t:= {d A_t\over dt},~ \b^\dbP_t := \sqrt{d\la M\ra_t\over dt},
\ea\right. \dbQ\mbox{-a.s.}
\eea
Similarly, for any $t\in [0, T)$, we may define $\cP_L^t$ on $\O^t$.

As in Denis, Hu and Peng \cite{DHP}, the set $\cP^t_L$ induces the following capacity:
 \bea
 \label{cC}
 \cC^L_t[A]
 :=
 \sup_{\dbP\in\cP^t_L}\dbP[A],
 ~~
 &\mbox{for all}&
 A\in\cF^t_T.
 \eea
We denote by $\dbL^1(\cF^t_T,\cP^t_L)$ the set of $\xi\in \dbL^0(\cF^t_T)$ satisfying $\sup_{\dbP\in\cP^t_L}\dbE^{\dbP}[|\xi|]<\infty$. The following nonlinear expectation will play a crucial role:
 \bea\label{cE}
 \overline{\cE}^L_t[\xi]
:=
 \sup_{\dbP\in\cP^t_L}\dbE^{\dbP}[\xi]
 ~~\mbox{and}~~
 \underline{\cE}^L_t[\xi]
 :=
 \inf_{\dbP\in\cP^t_L}\dbE^{\dbP}[\xi]
 =
 -\overline{\cE}^L_t[-\xi]
 ~~\mbox{for all}~~
 \xi\in\dbL^1(\cF^t_T,\cP^t_L).
 \eea
We remark that $\overline\cE^L[\xi]$ can be viewed as the solution of a Second Order BSDE (2BSDE, for short)  in the sense of \cite{STZ-2BSDE}, or a conditional $G$-expectation in the sense of \cite{Peng-G}. See Section \ref{sect-Example} for more details.  The following result will be important for us.

\begin{lem}
\label{lem-che}
For any $\ch\in \cH$ and any $L>0$, we have $\underline \cE^L_0[\ch] >0$.
\end{lem}
\proof By \reff{che}, we may assume $\ch_\e \le \ch$ for some $\e>0$. For any $\dbP\in \cP_L$ and $0<\d\le \e$, we have
\bea
\label{chest}
\dbP(\ch \le \d) \le \dbP(\ch_\e\le \d) = \dbP(\|B\|_\d \ge \e) \le \e^{-4} \dbE^{\dbP}[\|B\|_\d^4] \le CL^4 \e^{-4} \d^2.
\eea
This implies that, for $\d := {\e^2 \over \sqrt{2C}L^2}\wedge \e$,
\beaa
\dbE^\dbP[\ch] \ge \d \dbP(\ch > \d)  =\d\Big(1-\dbP(\ch \le \d) \Big) \ge \d \Big(1-CL^4 \e^{-4} \d^2) \ge {\d\over 2}.
\eeaa
Thus $\underline \cE^L_0[\ch] \ge {\d\over 2} >0$.
\qed

\begin{defn}
Let $X\in \dbL^0(\L)$ satisfy $X_t\in\dbL^1(\cF_t,\cP_L)$ for all $0\le t\le T$. We say that $X$ is an $\overline{\cE}^L-$supermartingale (resp. submartingale, martingale) if, for any $(t,\o)\in\L$ and any $\t\in \cT^t$,  $\overline{\cE}^L_t[X^{t,\o}_\t]\le$ (resp. $\ge,=$) $X_t(\o)$.
\end{defn}

We now state an important result for our subsequent analysis. Given a bounded process $X\in \dbL^0(\L)$, 
consider the nonlinear optimal stopping problem 
 \bea\label{cS}
 \overline{\cS}^L_t[X](\o) 
 := 
 \sup_{\t\in \cT^t}
 \overline\cE^L_t\big[X^{t,\o}_{\t}\big]
 &\mbox{and}&
 \underline{\cS}^L_t[X](\o) 
 := 
 \inf_{\t\in \cT^t}
 \underline\cE^L_t\big[X^{t,\o}_{\t}\big], 
 ~~(t,\o)\in\L.
 \eea
By definition, we have $\overline{\cS}^L[X]\ge X$ and $\overline{\cS}^L_T[X] = X_T$. The following nonlinear Snell envelope characterization is proved in \cite{ETZ0}.

\begin{thm}\label{thm-optimal}
Let $X\in\Usub$ be bounded, $\ch\in \cH$, and set $\wh X_t:=X_t\1_{\{t<\ch\}}+X_{\ch-}\1_{\{t\ge \ch\}}$. Define
 \beaa
 Y
 := 
 \overline\cS^L\big[\wh X\big]   
 &\mbox{and}& 
 \t^*:=\inf\{t\ge 0:Y_t=\wh X_t\}.
 \eeaa
Then $Y_{\t^*} = \wh X_{\t^*}$, $Y$ is an $\overline\cE^L$-supermartingale on $[0, \ch]$, and an $\overline\cE^L$-martingale on $[0, \t^*]$. Consequently, $\t^*$ is an optimal stopping time. 
\end{thm}

\begin{rem}\label{rem-h}
{\rm (i)
We emphasize that the maturity of the above nonlinear optimal stopping problem is restricted to be a hitting time in $\cH$. This requirement is due to technical aspects in the proof of Theorem \ref{thm-optimal} reported in \cite{ETZ0}. The difficulty is related to the regularity of the Snell envelope $Y$ and to some limiting arguments  under nonlinear expectation.

(ii) $Y$ is continuous  in $[0, \ch)$ and has a left limit at $\ch$.  However, in general $Y$ may have   a negative jump at $\ch$.}
\end{rem}

\subsection{The derivatives}

We define the path derivatives via the functional It\^{o} formula, which is initiated by  Dupire \cite{Dupire} and plays an important role in our paper.  
Denote 
\beaa
 \cP^t_\infty \;:=\; \bigcup_{L>0} \cP^t_L,
 ~~t\in[0,T].
\eeaa

\begin{defn}
\label{defn-spaceC12}  We say $u\in C^{1,2}(\L)$ if $u\in C^0(\L)$ and there exist  $\pa_t u \in C^0(\L)$, $\pa_\o u \in C^0(\L, \dbR^d)$, $\pa^2_{\o\o} u\in C^0(\L, \dbS^d)$ such that, for any 
$\dbP\in \cP^0_\infty$, $u$ is a $\dbP$-semimartingale satisfying:
\bea
\label{Ito}
d u = \pa_t u dt+ \pa_\o u \cd d B_t + \frac12 \pa^2_{\o\o} u : d \la B\ra_t,~~0\le t\le T,~~\dbP\mbox{-a.s.}
\eea
\end{defn}
We remark that the above  $\pa_t u$, $\pa_\o u$ and $\pa^2_{\o\o} u$, if they exist, are unique. Indeed, first considering $\dbP^{0,0}$, the probability measure corresponding to $\a = {\bf 0}, \b = {\bf 0}$ in \reff{cPL}, together with the required regularity  $\pa_t u \in C^0(\L)$ we obtain 
\bea
\label{pat}
\pa_t u(t,\o)  =  \lim_{h\downarrow 0} {1\over h}[u\big(t\!+\!h, \o_{\cd\wedge t}\big)       - u\big(t, \o\big)].
\eea
Next, considering $\dbP$ such that $\a^\dbP=1$, $\b^\dbP=0$ we obtain the uniqueness of $\pa_\o u$. Finally, considering $\dbP = \dbP_0$ we see that $\pa^2_{\o\o}u$ is also unique. Consequently, we call them the time derivative,  first order and second order space derivatives of $u$, respectively.  We define $C^{1,2}(\L^t)$   similarly. It is clear that, for any $(t,\o)$ and $u\in C^{1,2}(\L)$, we have $u^{t,\o}\in C^{1,2}(\L^t)$, and $\pa_t (u^{t,\o}) = (\pa_t u)^{t,\o}$, $\pa_\o (u^{t,\o}) = (\pa_\o u)^{t,\o}$, $\pa^2_{\o\o} (u^{t,\o}) = (\pa^2_{\o\o} u)^{t,\o}$.

\begin{rem}
\label{rem-pax}{\rm
(i) In Markovian case, namely $u(t,\o) = v(t,\o_t)$, if $v\in C^{1,2}([0, T]\times \dbR^d)$, then by the standard It\^{o} formula we see immediately that $u\in C^{1,2}(\L)$ with
\beaa
\pa_t u(t,\o) = \pa_t v(t,\o_t),\q \pa_\o u(t,\o) = D v(t,\o_t),\q \pa_{\o\o}^2 u(t,\o) = D^2 v(t,\o_t).
\eeaa
Here $D$ and $D^2$ denote the standard gradient and hessian of $v$ with respect to $x$. 

\no(ii) The typical case that the path derivatives exist is  that $u$ is smooth in Dupire's sense \cite{Dupire} (more precisely, the space $\dbC^{1,2}_b$ in Cont and Fournie \cite{CF}), and in that case our time derivative and space derivatives agree with the horizontal and vertical derivatives introduced therein, respectively, due to their functional It\^{o} formula. Therefore, any smooth function in the sense of Duprie's calculus is also smooth in the sense of  Definition \ref{defn-spaceC12}, namely our space $C^{1,2}(\L)$ is a priori larger than the  space $\dbC^{1,2}_b$ in \cite{CF}. In particular, Definition \ref{defn-spaceC12} is different from the corresponding definition in our previous paper \cite{EKTZ}, which uses Dupire's derivatives.

\no (iii) The main advantage of our definition is that all derivatives are defined within the continuous path space $\O$. In Dupire \cite{Dupire}, one has to extend the process $u$ (and the generator as well as the terminal condition of our PPDE later) to a larger domain $[0,T]\times \dbD([0,T])$, where $\dbD$ is the set of {\cad} paths. This is not natural in many situations, and is not necessary for our purpose, as it turns out that what we need is exactly the functional It\^{o} formula, rather than the precise form  of the derivatives. 

\no (iv) Moreover, compared to \cite{Dupire},   our definition does not require all the derivatives to be bounded, and we do not need \reff{Ito} to hold true for all semimartinagle measures $\dbP$. However, under our definition we do not require $\pa^2_{\o\o} u = \pa_\o (\pa_\o u) $. When $\pa_\o u$ is indeed differentiable, typically we should have $\pa^2_{\o\o} u = {1\over 2}\big[\pa_\o(\pa_\o u) + [\pa_\o(\pa_\o u)]^T\big]$.  For the last point see more details in \cite{BMZ}. 

\no (v) As explained in Cont and Fournie \cite{CF}, when $u$ is smooth enough in both senses, it holds that $\pa_\o u(t,\o) = D_t u(t,\o)$, where $D_t$ denotes the Malliavin derivative. We emphasize that, unlike the Malliavian derivative $D_t \xi$ which involves the perturbation of $\xi$ over the whole path $\o$,  $\pa_\o u(t,\o)$ involves the perturbation of $u$ only at the current time $t$. In particular, $\pa_\o u$ is $\dbF$-adapted. 
\qed}
\end{rem}

\begin{rem}
\label{rem-pax2}{\rm
 We shall also remark that, in our proof of comparison principle for PPDEs  in our accompanying paper \cite{ETZ2}, we actually uses only piecewise Markovian test functions and thus the standard It\^{o} formula is sufficient. So technically speaking, we can prove both the existence and uniqueness of viscosity solutions without using the path derivatives and the functional It\^{o} formula. However, it is more natural to consider truly path dependent test functions in this framework. In particular, it is more natural to define classical solutions for PPDEs by using path derivatives.  
\qed}
\end{rem}

\begin{eg}\label{eg-maximum}{\rm
Let $d=1$. As highlighted by Cont and Fournie \cite{CF}, a simple example of non-differentiable process is the running maximum process: $u(t,\o):=\overline{\o}_t:=\max_{0\le s\le t} \o_s$, $(t,\o)\in\L$.  Indeed, if it is differentiable, by \reff{pat} it is obvious that $\pa_t u(t,\o)=0$ for all $(t,\o)\in\L$. 
Then by \reff{Ito} one must have $\pa_\o u = 0$, and $\frac12\pa_{\o\o}^2 u dt = d \overline{B}_t$, which is impossible under $\dbP_0$. In terms of the Dupire's vertical derivatives, we have $\pa_\o u(t,\o)= 0$ whenever $\o_t<\overline{\o}_t$, and
 \beaa
 \pa_\o^+u(t,\o)=1
 &\mbox{and}&
 \pa_\o^-u(t,\o)=0
 ~~\mbox{whenever}~~
 \o_t=\overline{\o}_t,
 \eeaa
where $\pa_\o^+$ and $\pa_\o^-$ denote the right and left space derivatives in the sense of Dupire. Hence the process $u$ is not differentiable on $\{\o_t=\overline{\o}_t\}$.
\qed}
\end{eg}

\section{Fully nonlinear path dependent PDEs}
\label{sect-PPDE}
\setcounter{equation}{0}

In this paper we study the following fully nonlinear  parabolic path-dependent partial differential equation (PPDE, for short):
 \bea\label{PPDE}
 \cL u (t,\o)
 := \{-\pa_t u
          - G(., u, \pa_\o u, \pa^2_{\o\o} u)
    \}(t,\o) 
 =0,
 ~~0\le t< T,~\o\in \O,
\eea
where the generator $G :  \L \times  \dbR \times \dbR^d \times \dbS^d \rightarrow \dbR$ satisfies the following standing assumptions:

\begin{assum}\label{assum-G}
The nonlinearity $G$ satisfies:
\\
{\rm (i)} For fixed $(y,z,\g)$, $G(\cd, y,z,\g)\in \dbL^0(\L)$ 
and $|G(\cd, 0, {\bf 0}, {\bf 0})|\le C_0$. 
\\
{\rm (ii)} $G$ is elliptic, i.e. nondecreasing in $\g$. 
\\
{\rm (iii)} $G$ is  uniformly Lipschitz continuous in $(y,z,\g)$, with a Lipschitz constant $L_0$.
\\
{\rm (iv)} For any $(y, z,\g)$, $G(\cd, y,z,\g)$ is right continuous in $(t,\o)$ under $\dbf_\infty$, in the sense of Definition \ref{defn-rightcont}. 
\end{assum}

\begin{rem}
\label{rem-markovian}
{\rm In the Markovian case, namely $G(t, \o,.)= g(t, \o_t,.)$ and $u(t,\o) = v(t,\o_t)$, the PPDE \reff{PPDE} reduces to the following PDE: recalling Remark \ref{rem-pax} (i),
 \bea\label{PDE}
 \mathbf{L} v(t,x)
 := \{-\pa_t v
      -g(., v, D v, D^2 v)\}(t,x) 
 =0,~~ 
 t\in[0,T),~x\in \dbR^d.
 \eea
Namely, $u$ is a solution (classical or viscosity as we will introduce soon) of PPDE \reff{PPDE} corresponds to that  $v$ is a solution of PDE \reff{PDE}. However, slightly different from the PDE literature but consistent with \reff{pat}, here we should interpret $\pa_t v$ as  the right derivative of the function $v$ in $t$.
\qed}
\end{rem}

\subsection{Classical solutions}

\begin{defn}
\label{defn-classical} Let $u\in C^{1,2}(\L)$. We say $u$ is a classical solution   (resp. sub-solution, super-solution) of PPDE (\ref{PPDE}) if $\cL u (t,\o) =  \;(\mbox{resp.}\le,  \ge)\; 0$ for all $(t,\o) \in [0,T)\times \O$.
\end{defn}

It is clear that, in the Markovian setting as in Remark \ref{rem-markovian} with smooth $v$,
 $u$ is a classical solution   (resp. sub-solution,  super-solution) of PPDE (\ref{PPDE})
if and only if $v$ is a classical solution   (resp. sub-solution, super-solution) of PDE (\ref{PDE}).

\begin{eg}
\label{eg-classical1}{\rm 
 Let  $d=1$ and $u(t,\o):=\dbE^{\dbP_0}_t\big[\int_0^T B_tdt\big] (\o)= \int_0^t \o_sds + (T-t)\o_t$, $(t,\o)\in\L$. Then $u\in C^{1,2}(\L)$, and is a classical solution of the path dependent heat equation 
 \bea
 \label{heat}
 -\pa_tu-\frac12\pa^2_{\o\o}u=0
 \eea
  with terminal condition $u(T,\o)=\int_0^T \o_tdt$.
\qed}
\end{eg}

\begin{eg}
\label{eg-classical2}{\rm 
Let $d=1$ and $u(t,\o):=\dbE^{\dbP^t_0}\Big[\overline{(\o\otimes_tB^t)}_T\Big]$, $(t,\o)\in\L$ with the notation of Example \ref{eg-maximum}. Then one can easily check that $u(t,\o) = v(t,\o_t, \overline \o_t)$, where $v$ is a deterministic function defined by:
 \bea
 \label{v}
 \ba{rcl}
 v(t,x,y)
 &:=& 
 \dbE^{\dbP^t_0} \Big[ y \vee (x  + (\overline {B^t})_T)\Big] 
 \;=\; 
 x + \sqrt{T-t} \psi(\frac{y-x}{\sqrt{T-t}}),\q x\le y
 \\
 \psi(z) &:=& 
 \dbE^{\dbP_0} \Big[ z \vee \overline B_1\Big] 
 \;=\;  \dbE^{\dbP_0} \Big[ z \vee  |B_1|\Big] 
 \;=\;  z[2\Phi(z) -1]+ \frac{2}{\sqrt{2\pi}}e^{-z^2/2},
 ~~ z\ge 0,
 \ea
 \eea
and $\Phi$ denotes the cdf of  the standard normal distribution.  
We note that  $v$ is smooth for $t<T$, and $D_y v (t,x,x) =0$. Since the support of $d \overline B_t$ is in $\{B_t = \overline B_t\}$, it follows that $D_y v(t, B_t, \overline B_t) d\overline{B}_t=0$. This implies that
\beaa
d u(t,\o) = d v(t, \o_t, \overline \o_t) 
          = \pa_t v dt + D_x  v dB_t + \frac12 D^2_{xx}  v d\la B\ra_t.
\eeaa
By \reff{pat} it is clear that $\pa_t u(t,\o) = \pa_t v(t, \o_t, \overline \o_t)$. Then by \reff{Ito} we see that $\pa_\o u(t,\o) = D_x  v(t, \o_t, \overline \o_t)$ and $\pa_{\o\o}^2 u(t,\o) = D^2_{xx}  v(t, \o_t, \overline \o_t)$. Thus $u \in C^{1,2}(\L)$.

Finally, it is straightforward to check that $u$ is a classical solution to the path dependent heat equation   \reff{heat} with terminal condition $u(T,\o)=\overline B_T$.
\qed}
\end{eg}

\begin{rem}
\label{rem-classical}{\rm We shall remark that, unlike a standard heat equation which always has classical solution in $[0, T)$, a path dependent one may not have a classical solution in $[0, T)$. One simple example is the equation \reff{heat} with terminal condition $u(T,\o) = B_{t_0}(\o)$ for some $0< t_0<T$. Then clearly $u(t, \o) = B_{t\wedge t_0}(\o)$, and thus $\pa_\o u(t, \o) = \1_{[0, t_0]}(t)$ is discontinuous. Following Proposition \ref{prop-semi} below and our accompanying paper \cite{ETZ2} (or Section \ref{sect-semilinear} below under a slight reformulation), and weakening the boundedness assumption as pointed out in Remark \ref{rem-bound} below, one can easily see $u$ is the unique viscosity solution of the  equation \reff{heat} with terminal condition $u(T,\o) = B_{t_0}(\o)$.  We refer to Peng and Wang \cite{PW} for  sufficient conditions of existence of classical solutions for more general semilinear PPDEs.
\qed}
\end{rem}

\subsection{Definition of viscosity solutions}
\label{sect-defvisco}

We next introduce our notion of viscosity solutions. Recall the nonlinear Snell envelope notation \reff{cS}. For any  $u\in \dbL^0(\L)$, $(t,\o) \in [0, T)\times \O$, and  $L>0$,
define
 \begin{equation}\label{cA}
 \ba{lll}
\!\underline\cA^{\!L}\!u(t,\o) 
\! :=\!
 \Big\{\f\in C^{1,2}(\L^{\!t}):
       (\f-u^{t,\o})_t 
       = 0 =
      \underline\cS^L_t\big[(\f-u^{t,\o})_{\cdot\wedge\ch}
                       \big]
      ~\mbox{for some}~\ch\in \cH^t
 \Big\},
 \\
\! \overline\cA^{\!L}\!u(t,\o) 
 \!:=\! 
 \Big\{\f \in C^{1,2}(\L^{\!t}):
      (\f-u^{t,\o})_t 
      =0=
      \overline\cS^L_t\big[(\f-u^{t,\o})_{\cdot\wedge\ch}
                       \big] 
      ~\mbox{for some}~\ch\in \cH^t
 \Big\}.
 \ea
 \end{equation}

\begin{defn}
\label{defn-viscosity}
\no {\rm (i)} Let $L>0$. We say $u\in\Usub$ (resp. $\Usup$) is a viscosity $L$-subsolution (resp. $L$-supersolution) of PPDE (\ref{PPDE})  if,  for any $(t,\o)\in [0, T)\times \O$ and any $\f \in \underline\cA^{L}u(t,\o)$ (resp. $\f \in \overline\cA^{L}u(t,\o)$):
 \beaa
 \cL^{t,\o}\f(t,{\bf 0}) :=  \big\{-\pa_t \f   -G^{t,\o}(., \f,\pa_\o \f,\pa^2_{\o\o}\f)  
 \big\}(t,{\bf 0})
 &\le  ~~(\mbox{resp.} \ge)&  0.
 \eeaa

\no {\rm (ii)} $u\in\Usub$ (resp. $\Usup$)  is a viscosity subsolution (resp. supersolution) of PPDE (\ref{PPDE}) if  $u$ is viscosity $L$-subsolution (resp. $L$-supersolution) of PPDE (\ref{PPDE}) for some $L>0$.

\no {\rm (iii)} $u\!\in\! UC_b(\L)$ is viscosity solution of PPDE (\ref{PPDE}) if it is viscosity sub- and supersolution.
\end{defn}

\begin{rem}
\label{rem-bound}
{\rm For technical simplification, in this paper and the accompanying one \cite{ETZ2}, we consider only bounded viscosity solutions. By some more involved estimates one can extend our theory to viscosity solutions satisfying certain growth conditions. We shall leave this for future research, however, in some examples below we may consider unbounded viscosity solutions as well.
\qed}
\end{rem}

\begin{rem}
\label{rem-right}
{\rm Since our PPDE is backward, in \reff{cA} the test functions $\f$ are defined only after $t$. By this nature, both the viscosity solution $u$ and the generator $G$ are required only to be right continuous in $(t,\o)$ under $\dbf_\infty$. To prove the comparison principle, however, we will assume some stronger regularity of $G$, see our accompanying paper \cite{ETZ2}.
\qed}
\end{rem}
We next provide an intuitive justification of our Definition \ref{defn-viscosity} which shows how the above nonlinear optimal stopping problems $\underline{\cS}$ and $\overline{\cS}$ appear naturally. 

Let $u\in C^{1,2}(\L)$ be a classical supersolution of PPDE (\ref{PPDE}), $(t^*,\o^*) \in [0, T)\times \O$, and $\f\in C^{1,2}(\L^{t^*})$. Then:
 \bea
 0
 \le \cL u(t^*,\o^*) = \cL^{t^*,\o^*}\f(t^*,{\bf 0})+R(t^*,{\bf 0})
 \label{consistency-intuition}
 \eea
where $R(t,\o)= \pa_t(\f-u^{t^*, \o^*})(t,\o)
 +\hat\a\cdot\pa_\o(\f-u^{t^*, \o^*})(t,\o)
 +\frac12\hat\b^2:\pa^2_{\o\o}(\f-u^{t^*, \o^*})(t,\o)$ for $(t,\o)\in \L^{t^*}$, $\hat\alpha:= G_z(t^*,\o^*,u(t^*,\o^*),\hat z,\hat\g),$ and $\hat\b:=\big(2G_\g(t^*,\o^*,u(t^*,\o^*),\hat z,\hat\g)\big)^{1/2}$ are constant drift and diffusion coefficients, and $(\hat z,\hat\g)$ are some convex combination of $(\partial_\o u,\partial^2_{\o\o}u)(t^*,\o^*)$ and $(\partial_\o\varphi,\partial^2_{\o\o}\varphi)(t^*,{\bf 0})$. 
 
The question is how to choose the test process $\varphi$ so as to deduce from \reff{consistency-intuition} that $ \cL^{t^*,\o^*}\f(t^*,0) \ge 0$. A natural sufficient condition is $R(t^*,{\bf 0}) \ge 0$. To achieve that, 
our crucial observation is that 
 \beaa
 d(\f-u^{t^*,\o^*})(t,\o)
 &=&
 R(t,\o)dt+\pa_\o(\f-u^{t^*,\o^*})(t,\o)\cdot\hat\b d\hat W_t,
 ~~\hat\dbP-\mbox{a.s.}
 \eeaa
where $\hat W$ is a Brownian motion under the probability measure $\hat\dbP\in\cP^{t^*}_{L_0}$ defined by the pair $(\hat\alpha,\hat\beta)$, and $L_0$ is the Lipschitz constant of the nonlinearity $G$. Therefore, in order to ensure  $R(t^*, {\bf 0})\le 0$, we have to choose the test process $\f$ so that the difference $(\f-u^{t^*,\o^*})$ has a nonpositive $\hat\dbP-$drift locally at the right hand-side of $t^*$. This essentially means that $(\f-u^{t^*,\o^*})$ is a $\hat\dbP-$supermartingale on some right-neighborhood $[t^*,\ch]$ of $t^*$, and therefore $(\f-u^{t^*,\o^*})_{t^*}\ge\dbE^{\hat\dbP}[(\f-u^{t^*,\o^*})_{\t\wedge\ch}]$ for any stopping time $\t$. Since the probability measure $\hat\dbP$ is imposed by the above calculation, we must choose the test process $\f$ so that
$(\f-u^{t^*,\o^*})_{t^*}\ge
 \overline{\cE}^L_{t^*}
 \big[(\f-u^{t^*,\o^*})_{\t\wedge\ch}\big]$ for all stopping time $\tau$. Finally, since $\tau=t^*$ is a legitimate stopping rule, we arrive at 
 \beaa
 (\f-u^{t^*,\o^*})_{t^*}
 &=&
 \overline{\cS}^L_{t^*}
 \big[(\f-u^{t^*,\o^*})_{\t\wedge\ch}\big],
 \eeaa
which corresponds exactly to our definition of $\overline{\cA}^Lu(t^*,\o^*)$.  

Conversely, if the pair $\big((t^*,\o^*),\varphi\big)$ satisfies the last equality, then it follows from the Snell envelope characterization of Theorem \ref{thm-optimal} that $(\f-u^{t^*,\o^*})_{t^*}\ge\overline{\cE}^L_{t^*}\big[\overline{\cS}^L_{\tau\wedge\ch}[\f-u^{t^*,\o^*}]\big]\ge\overline{\cE}^L_{t^*}\big[(\f-u^{t^*,\o^*})_{\tau\wedge\ch}\big]$, for all stopping time $\tau$. By the right-continuity, this implies that  $R(t^*, {\bf 0})\le 0$. Hence our definition of the set of test processes $\overline{\cA}^Lu(t^*,\o^*)$ is essentially necessary and sufficient for the inequality $R(t^*,{\bf 0})\le 0$.

\begin{rem}
\label{rem-semilinear}
{\rm From the last intuitive justification of our definition, we see that for a semilinear path-dependent PDE, $\hat\beta$ is a constant matrix. Then, in agreement with our previous paper \cite{EKTZ}, also see Section \ref{sect-semilinear} below,  it is not necessary to vary the coefficient $\beta$ in the definition of the operator $\overline{\cE}^L$.

Similarly, in the context of a linear PPDE, both coefficients $\hat\alpha$ and $\hat\beta$ are constant, and we may define the sets $\overline{\cA}^Lu$ and $\underline{\cA}^Lu$ by means of the linear expectation operator. Finally, for a first order PPDE, we may take the diffusion coefficient $\beta\equiv 0$, see Section \ref{sect-First}.
\qed}
\end{rem}

In the rest of this section we provide several remarks concerning our definition of viscosity solutions. In most places we will comment on the viscosity subsolution only, but obviously similar properties hold for the viscosity supersolution as well.

\begin{rem}
\label{rem-cA}
{\rm As standard in the literature on viscosity solutions of  PDEs:
\\
(i) The viscosity property is a local property in the following sense. For any $(t,\o) \in [0,T)\times \O$ and any $\e>0$, define as in \reff{che},
\bea
\label{chet}
\ch^t_\e := \inf\Big\{ s>t: |B^t_s| \ge \e\Big\} \wedge (t+\e) &\mbox{and thus}& \ch_\e = \ch^0_\e.
\eea
It is clear that $\ch^t_\e \in \cH^t$.   To check the viscosity property of $u$ at $(t,\o)$, it suffices to know the value of $u^{t,\o}$ on $[t, \ch_\e]$ for an arbitrarily small $\e>0$. In particular, since $u$ and $\f$ are locally bounded, there is no integrability issue in \reff{cA}. Moreover,  for any $\f\in \underline \cA^L u(t,\o)$ with corresponding $\ch\in \cH^t$,  by \reff{che} we have $\ch^t_\e \le \ch$ when $\e$ is small enough.
\\
(ii)  The fact that $u$ is a viscosity solution does not mean that the PPDE must hold with equality at some $(t,\o)$ and $\f$ in some appropriate set. One has to check viscosity subsolution property and  viscosity supersolution property separately.
\\
(iii) In general  $\underline\cA^{L}u(t,\o)$ could be empty. In this case automatically $u$ satisfies the viscosity subsolution property at $(t,\o)$. 
\qed}
\end{rem}

\begin{rem}
\label{rem-markovian-viscosity}
{\rm (i) Consider the Markovian setting in Remark \ref{rem-markovian}. One can easily check that $u$ is a viscosity subsolution of PPDE \reff{PPDE} in the sense of Definition \ref{defn-viscosity} implies that $v$ is a viscosity subsolution of PDE \reff{PDE} in the standard sense, see e.g. \cite{CIL} or \cite{FS}. However, the opposite direction is in general not true. We shall point out though,  when the PDE is wellposed, by uniqueness our definition of viscosity solution of PPDE \reff{PPDE}  is consistent with the  viscosity solution of PDE \reff{PDE} in the standard sense.  Moreover, we emphasize that our definition involves a richer set of test functions which in principle opens the hope for an easier proof of uniqueness.
\\
(ii) Definition \ref{defn-viscosity} does not reduce to the definition introduced in the semilinear context of \cite{EKTZ} (or Section \ref{sect-semilinear} below) either, because we are using a different nonlinear expectation $\overline\cE^L$ here. It is obvious that any viscosity subsolution in the sense of \cite{EKTZ} is also a viscosity subsolution in the sense of this paper, but the opposite direction is in general not true. However, the definitions of viscosity solutions  are actually equivalent for semilinear PPDEs,  in view of the uniqueness result of our accompanying paper \cite{ETZ2}. See also Remark \ref{rem-semilinear}.
\qed}
\end{rem}

\begin{rem}
\label{rem-L}
{\rm For $0< L_1 < L_2$,  obviously  $\cP_{L_1}^t \subset \cP_{L_2}^t$, $\underline\cE^{L_2}_t\le \underline\cE^{L_1}_t$, and $\underline\cA^{L_2}u(t,\o) \subset \underline\cA^{L_1}u(t,\o)$.  Then one can easily check that a viscosity $L_1$-subsolution must be a viscosity $L_2$-subsolution. Consequently, $u$ is a viscosity subsolution if and only if

\ms
\q there exists an $L\ge 1$ such that, for all $L'\ge L$, $u$ is a viscosity $L'$-subsolution.
{ \hfill \vrule width.25cm height.25cm depth0cm\smallskip}}
\end{rem}

We next report the following result whose proof follows exactly the lines of Remark 3.9 (i) in \cite{EKTZ}.

\begin{prop}\label{prop-change} 
Let Assumption \ref{assum-G} hold true, and let $u$ be a viscosity subsolution of PPDE \reff{PPDE}. For $\l\in \dbR$, the process $\tilde u_t := e^{\l t}u_t$ is a viscosity subsolution of:
\bea
\label{tildePPDE}
\tilde\cL \tilde u
&:=&
- \pa_t \tilde u  -\tilde G(t, \o, \tilde u, \pa_\o\tilde u, \pa^2_{\o\o} \tilde u)
\;\le\; 0,
\eea
where $\tilde G(t, \o, y,z,\g):=-\l y+e^{\l t} G(t,\o, e^{-\l t}y, e^{-\l t} z, e^{-\l t} \g)$.
\end{prop}

\begin{rem}\label{rem-chgevariable}{\rm
Under Assumption \ref{assum-G}, we are not able to prove a more general change of variable formula. However, this will be achieved under stronger assumptions, see Proposition 4.5 and Theorem 4.6 of our accompanying paper \cite{ETZ2}.
\qed}
\end{rem}

\subsection{Consistency with classical solutions}
\label{sect-consistency}

\begin{thm}
\label{thm-consistency}
Let Assumption \ref{assum-G} hold and $u \in C^{1,2}(\L)\cap UC_b(\L)$. Then $u$ is a classical solution (resp. subsolution, supersolution)  of PPDE \reff{PPDE} if and only if it is a viscosity solution (resp.  subsolution, supersolution).
\end{thm} 
\proof We prove the subsolution property only. 
Assume $u$ is a viscosity $L$-subsolution. For any $(t,\o)$,  since $u \in C^{1,2}(\L)$, we have $u^{t,\o}\in C^{1,2}(\L^t)$ and thus $u^{t,\o} \in \underline \cA^Lu(t,\o)$ with $\ch := T$. By definition of viscosity $L$-subsolution we see that $\cL u(t,\o) \le 0$. 

On the other hand, assume $u$ is a classical subsolution. If $u$ is not a viscosity subsolution, then it is not a viscosity $L_0$-subsolution. Thus there exist $(t,\o)\in\L$ and $\f\in {\underline \cA}^{L_0}u(t,\o)$ 
such that $2c := \cL \f(t, {\bf 0}) >0$. Without loss of generality, we set $t:=0$ and, by Remark \ref{rem-cA} (i),  let $\ch = \ch_\e \in \cH$ defined in \reff{che} for some small constant $\e>0$ be the hitting time used  in the definition of  ${\underline \cA}^{L_0}u(0,{\bf 0})$. Now recall \reff{cPL} and let $\dbP\in \cP_{L_0}$ corresponding to some constants $\a\in \dbR^d$ and $\b\in\dbS^d$ which will be determined later. Then
 \beaa
 0 
 &\le& 
 \underline\cE^{L_0}_0\big[(\f- u)_{\ch} \big] 
 \le  
 \dbE^\dbP\big[(\f- u)_{\ch} \big].
\eeaa
Applying functional It\^o's formula \reff{Ito} and  noticing that $(\f- u)_0=0$, we have
 $$
 (\f- u)_{\ch} 
 = 
 \!\!\int_0^{\ch}\!\! \Big[\pa_t(\f-u)_s 
                  + \frac12\pa^2_{\o\o}(\f-u)_s\!:\!\b^2  
                  + \pa_\o(\f-u)_s \cd \a\Big] ds 
            + \!\!\int_0^{\ch}\!\!\pa_\o(\f-u)_s\!\cd\!\b dW^\dbP_s.
 $$
Taking expected values, this leads to
 $$
 0
 \le 
 \dbE^\dbP\Big[ \int_0^{\ch} 
                \Big(\pa_t(\f-u)_s 
                     + \frac12\pa^2_{\o\o}(\f-u)_s:\b^2
                     + \pa_\o(\f-u)_s\cd\a
                \Big)ds
          \Big]
 =
 \dbE^\dbP\Big[ \int_0^{\ch} (\tilde \cL\f - \tilde \cL u)_s ds\Big],
 $$
where $
\tilde\cL \f_s := -\cL\f_s  - G(\cd, \f, \pa_\o \f, \pa^2_{\o\o}\f)_s + \frac12(\pa^2_{\o\o}\f)_s : \b^2  + (\pa_\o \f )_s \cd \a$. 
 Since $\tilde\cL \f$ and $\tilde\cL u$ are continuous, for $\e$ small enough we have 
$|\tilde \cL\f_s - \tilde \cL\f_0| + |\tilde \cL u_s - \tilde \cL u_0| \le \frac{c}{2}$ on $[0,\ch]$. Then
\bea
\label{consistency-est}
0&\le&  \dbE^\dbP\Big[ (\tilde \cL\f_0 - \tilde \cL u_0 + c)\ch\Big].
\eea
Note that $\cL u_0 \le 0$, $\cL\f_0 = 2c$, and $\f_0 = u_0$. Thus
\beaa
\tilde \cL\f_0 -  \tilde \cL u_0 
&\le& 
-2c + \frac12\pa^2_{\o\o}(\f-u)_0:\b^2
+ \pa_\o(\f-u)_0\cd \a \\
&&- [G\big(., u, \pa_\o \f, \pa^2_{\o\o}\f\big)_0  
- G\big(., u, \pa_\o u, \pa^2_{\o\o}u\big)_0 ].
\eeaa
By Assumption \ref{assum-G} (iii), there exist $\a$ and $\b$ such that $\dbP\in \cP_{L_0}$ and
$$
G\big(., u, \pa_\o \f, \pa^2_{\o\o}\f\big)_0  
- G\big(., u, \pa_\o u, \pa^2_{\o\o}u\big)_0
=
\frac12\pa^2_{\o\o}(\f-u)_0 : \b^2  
+ \pa_\o (\f-u)_0\cd\a.
$$
Then $\tilde \cL\f_0 - \tilde \cL u_0 \le -2c$, and \reff{consistency-est} leads to $0\le\dbE^\dbP[ -c \ch] < 0$, contradiction.
\qed

\section{Some Examples with (Semi-)explicit Solution}
\label{sect-Example}
\setcounter{equation}{0}

In this section, we study several special PPDEs which have (semi-)explicit viscosity solutions, for example via backward SDEs or second order BSDEs. These solutions provide probabilistic representations for the PPDEs and thus can be viewed as path dependent nonlinear Feynman-Kac formula. More importantly, as value functions of some stochastic control problems, these examples illustrate how to check the viscosity properties of processes arising in applied problems. As in the viscosity theory of PDEs, the main tools are the regularity of the processes in $(t,\o)$ and the dynamic programming principle.

\subsection{First order PPDEs}
\label{sect-eg-first}

\begin{eg}\label{eg-dtphi=0}
{\rm  Suppose that $u(t,\o)=v(\o_t)$ for all $(t,\o)\in\Lambda$, where $v : \dbR^d \to \dbR$ is bounded and continuous. Then by \reff{pat} we should have $\pa_tu=0$. We now verify that $u$ is a viscosity solution of the equation $-\partial_tu=0$.

Indeed, for $\varphi\in\overline{\cA}^Lu(t,\o)$, it follows from our definition that, for some $\ch\in \cH^t$:  
 $$
 (\varphi-u^{t,\o})_t = 0
 \ge
 \dbE^{\dbP^{0,0}}
 \big[(\varphi-u^{t,\o})_{(t+\d)\wedge \ch}\big] ~~\mbox{for all}~\d>0.
 $$
where $\dbP^{0,0}$ is  again  the probability measure corresponding to $\a = {\bf 0}, \b = {\bf 0}$ in \reff{cPL}. 
Notice that under $\dbP^{0,0}$, the canonical process $\o$ is frozen to its value at time $t$. Then $\ch=T$, $\dbP^{0,0}$-a.s. and thus, for $\d < T-t$, 
 \beaa
 \f(t, {\bf 0}) - v(\o_t) =   (\varphi-u^{t,\o})_t
 \ge
\dbE^{\dbP^{0,0}}
 \big[(\varphi-u^{t,\o})_{(t+\d)\wedge \ch}\big] = \f(t+\d, {\bf 0}) - v(\o_t).
\eeaa
 This implies that $\pa_t\varphi(t,{\bf 0})\le 0$. A similar argument shows that $\pa_t\f(t,{\bf 0})\ge 0$ for all $\f\in\underline{\cA}^Lu(t,\o)$.
\qed}
\end{eg}

\begin{eg}\label{eg-maximum2}{\rm 
Let $d=1$ and use the notations in Example \ref{eg-maximum}.  We  check that $u(t,\o) := 2 \overline B_t - B_t$ is a viscosity solution of the first order equation:
 \bea
 \label{maxPPDE}
 -\pa_t u -|\pa_\o u| + 1
 &=&
 0.
 \eea
 By Example \ref{eg-maximum}, $u$ is not smooth, so it is a viscosity solution but not a classical solution.
 
 When $\o_t < \overline\o_t$, it is clear that $u$ is smooth with $\pa_t u (t,\o) = 0, \pa_\o u(t, \o) = -1$ and thus satisfies \reff{maxPPDE}. So it suffices to check the viscosity property when $\o_t = \overline \o_t$. Without loss of generality, we check it at $(t,\o) = (0, 0)$. 
 \\
(i) We first check that ${\underline\cA}^L u(0, 0)$ 
 is empty for $L\ge 1$, and thus $u$ is a viscosity subsolution. Indeed, assume $\f\in {\underline\cA}^L u(0, 0)$ with corresponding $\ch\in \cH$. By Remark \ref{rem-cA} (i), without loss of generality we may assume $\ch = \ch_\e$ for some small $\e>0$, and thus $\pa_t \f, \pa^2_{\o\o}\f$ are bounded on $[0, \ch]$. Note that $\dbP_0\in \cP_L$. By definition of ${\underline\cA}^L$ we have, for any $0<\d<\e$,
\beaa
0 
&\le& 
\dbE^{\dbP_0}\Big[(\f - u)_{\d\wedge \ch}\Big] 
=  
\dbE^{\dbP_0}\Big[\int_0^{\d\wedge \ch}(\pa_t \f + \pa_{\o\o}^2\f)(t,\o)ds
                  - 2 \overline B_{\d\wedge \ch}
             \Big]
\\
&\le& 
C\dbE^{\dbP_0}[\d\wedge \ch] - 2  \dbE^{\dbP_0}[ \overline B_{\d\wedge \ch}] \le C \d - 2  \dbE^{\dbP_0}[ \overline B_{\d}] +  2\dbE^{\dbP_0}[ \overline B_\d \1_{\{\ch \le \d\}}] \\
&\le& C\d - c\sqrt{\d} + C\sqrt{ \dbP_0(\ch \le \d)} \le C[\d +  \e^{-2}\d] - c\sqrt{\d},
\eeaa
where $c := 2\dbE^{\dbP_0}[ \overline B_1]>0$ and the last inequality thanks to \reff{chest}.
This leads to a contradiction when $\d$ is small enough. Therefore, ${\underline\cA}^L u(0, {\bf 0})$ is empty.
\\
(ii) We next check the viscosity supersolution property.  Assume to the contrary that $-c :=-\pa_t \f(0,0)  -|\pa_\o \f(0, 0)| + 1 < 0$ for some $\f \in {\overline\cA}^L u(0, 0)$ and $L\ge 1$. Let $\a := \sgn (\pa_\o \f(0, 0))$ (with the convention $\sgn(0) := 1$), $\b := 0$, and $\dbP\in \cP_L$ be determined by \reff{cPL}. When $\a=1$, we have $B_t = t, \overline B_t = t$, $\dbP$-a.s. When $\a=-1$, we have $B_t = -t$, $\overline B_t = 0$, $\dbP$-a.s. In both cases, it holds that  $u(t,\o) = t$, $\ch_\e = \e$, $\dbP$-a.s. By choosing $\ch = \ch_\e$ and $\e$ small enough, we may assume $|\pa_t \f(t, B) - \pa_t \f(0,0)| + |\pa_\o\f (t, B) - \pa_\o\f(0,0)|  \le \frac{c}{2}$ for $t \le \ch_\e$.  By the definition of $ {\overline\cA}^L u(0, 0)$ we get
 \beaa
 0 
 &\ge& 
 \dbE^{\dbP}\Big[(\f - u)_{\ch_\e}\Big] 
 =  \dbE^{\dbP}\Big[ \int_0^\e (\pa_t \f + \a \pa_\o \f )_t dt 
                    - \e\Big] \\
 &\ge&
 \dbE^{\dbP}\Big[ \int_0^\e 
                  \Big(\pa_t\f_0+\a\pa_\o\f_0- \frac{c}{2}
                  \Big) dt 
            \Big] -\e
 \\
 &=& 
 \dbE^{\dbP}\Big[ \int_0^\e 
                  \Big(\pa_t\f_0+|\pa_\o\f_0| - \frac{c}{2}
                  \Big) dt 
 \Big] -\e 
 \;=\;  \int_0^\e \Big(1+c- \frac{c}{2}\Big) dt -\e 
 \;=\; \frac12\;c\e
 \;>\;
 0.
\eeaa
This is the required contradiction, and thus $u$ is a viscosity supersolution of \reff{maxPPDE}.
\qed
}
\end{eg}

\subsection{Semi-linear PPDEs and BSDEs}

We  now consider the following semi-linear PPDE:
\bea
\label{semiPPDE}
-\pa_t u - \frac{1}{2} \si^2(t,\o): \pa_{\o\o}^2 u 
- F\big(t,\o, u, \si(t,\o)\pa_\o u \big) =0,
\q u(T,\o) = \xi(\o),
\eea
where $\si \in \dbL^0(\dbF, \dbS^d)$, $\xi\in \dbL^0(\cF_T)$, and $F$ is $\dbF$-progressively measurable in all variables. We note that \cite{EKTZ} studied the case $\si =I_d$ for simplicity.   We shall assume

\begin{assum}
\label{assum-F}
{\rm (i)} $\si$, $F(t,\o, 0, {\bf 0})$, and $\xi$ are bounded by $C_0$, and $\si>{\bf 0}$.
\\
{\rm (ii)}  $\si$  is uniformly Lipschitz continuous in $\o$ and $F$ is uniformly Lipschitz contin. in $(y,z)$.
\\
{\rm (iii)} $F$ and $\xi$ are uniformly continuous in $\o$, and the common modulus of continuity function $\rho_0$ has polynomial growth. 
\\
{\rm (iv)} $\si$ and $F(\cd, y, z)$ are right continuous in $(t,\o)$ under $\dbf_\infty$ for any $(y,z)$, in the sense of Definition \ref{defn-rightcont}. 
\end{assum}

The boundedness in Assumption \ref{assum-F} (i) is just for simplification, and can be weakened to some growth condition. The assumption $\si > {\bf 0}$ and that $F$ depends on the gradient term through the special form $\si(t,\o)\pa_\o u$ are mainly needed for the subsequent BSDE representation. 
 
For any $(t,\o)\in\L$, consider the following decoupled FBSDE on $[t,T]$:
 \bea
 \label{BSDEstrong}
 \left\{\ba{lll}
 \dis \cX_s 
 = 
 \int_t^s \si^{t,\o}(r, \cX_\cd)d B^t_r,
 \\
 \dis  \cY_s 
 = 
 \xi^{t,\o}(\cX) + \int_s^T F^{t,\o}(r,\cX_\cd,\cY_r,\cZ_r)dr 
 - \int_s^T \cZ_r\cdot dB^t_r,
 \ea\right.
 ~\dbP^t_0-\mbox{a.s.}
 \eea
Under Assumption \ref{assum-F}, clearly FBSDE \reff{BSDEstrong} has a unique solution $(\cX^{t,\o}, \cY^{t,\o}, \cZ^{t,\o})$.  Alternatively, we may consider the BSDE in weak formulation:
 \begin{equation}
 \label{BSDEweak}
 Y^{t,\o}_s 
 = 
 \xi^{t,\o}(B^t) 
 + \int_s^T F^{t,\o}(r, B^t_\cd,Y^{t,\o}_r, Z^{t,\o}_r) dr 
 - \int_s^T Z^{t,\o}_r\cdot(\si^{t,\o}(r, B^t_\cd))^{-1}dB^t_r,
~\dbP^{t,\o}\mbox{-a.s.}
 \end{equation}
where $\dbP^{t,\o} := \dbP_0^t \circ(\cX^{t,\o})^{-1}$ denotes the distribution of  $\cX^{t,\o}$. Then, for any fixed $(t,\o)$,
\beaa
\cY^{t,\o}_t = Y^{t,\o}_t &\mbox{and}& \mbox{is a constant due to the Blumenthal zero-one law}.
\eeaa

\begin{prop}
\label{prop-semi}
Under Assumption   \ref{assum-F}, $u(t,\o) := Y^{t,\o}_t = \cY^{t,\o}_t$ is a viscosity solution of PPDE \reff{semiPPDE}.
\end{prop}

\proof  We proceed in two steps.
\\
{\it Step 1.} In Step 2 below, we will show that $u\in UC_b(\L)$ and satisfies the dynamic programming principle:  for any $(t,\o)\in \L$ and $\t\in\cT^t$,
 \bea\label{semiDPP1}
 Y^{t,\o}_s 
\! = \!
 u^{t,\o}(\t, B^t) 
 + \!\int_s^\t F^{t,\o}(r, B^t_\cd,  Y^{t,\o}_r, Z^{t,\o}_r) dr 
 -\! \int_s^\t Z^{t,\o}_r\cdot(\si^{t,\o}(r, B^t_\cd))^{-1}dB^t_r,
 \dbP^{t,\o}\mbox{-a.s.}
\eea
Let $L$ be a  Lipschitz constant of $F$ in $z$ satisfying $|\si| \le \sqrt{2L}$. We now show that  $u$ is an $L$-viscosity solution. Without loss of generality, we prove only the viscosity subsolution property at $(t,\o) = (0, {\bf 0})$.  For notational simplicity we omit the superscript $^{0,{\bf 0}}$ in the rest of this proof. Assume to the contrary that, 
  \beaa
 c:= 
 -\big\{\pa_t \f + \frac{1}{2} \si^2: \pa_{\o\o}^2 \f 
        +F(\cd, u, \si\pa_\o\f)\big\}(0,{\bf 0}) 
 > 0 
 &\mbox{for some}&\f\in \underline \cA^{L} u(0,{\bf 0}).
 \eeaa
Let $\ch\in \cH$ be the hitting time corresponding to $\f$ in  \reff{cA}, and by Remark \ref{rem-cA} (i), without loss of generality we may assume $\ch = \ch_\e$ for some small $\e>0$. Since $\f \in C^{1,2}(\L)$ and $u\in UC_b(\L)$, by Assumption \ref{assum-F} (iv) and the uniform Lipschitz property of $F$ in $(y,z)$, we may assume $\e$ is small enough such that 
 \beaa
 -\big\{\pa_t \f + \frac{1}{2} \si^2: \pa_{\o\o}^2 \f 
        + F(\cd, u, \si\pa_\o\f)\big\}(t,\o) 
 \ge \frac{c}{2} >0,\q t\in [0, \ch].
 \eeaa
 Notice that $d \la B \ra_t = \si^2(t, B_\cd)dt$, $\dbP$-a.s. Using the dynamic programming principle \reff{semiDPP1}, and applying It\^{o}'s formula on $\f$, we have:
 \beaa
 &&
 (\f - u)_\ch = (\f - u)_\ch -(\f - u)_0 \\
 &=& 
\int_0^\ch \big(\pa_\o \f -  \si^{-1}Z\big)(s,B_\cd)\cdot dB_s+ \int_0^\ch\Big(\pa_t\f + \frac12 \si^2 : \pa_{\o\o}^2\f 
                  + F(\cd, u, Z)\Big)(s, B_\cd) ds 
 \\
 &\le&
 \int_0^\ch \big(\si\pa_\o \f -Z\big)(s, B_\cd)\cdot \si^{-1}(s, B_\cd)dB_s- \int_0^\ch\Big( \frac{c}{2} 
                 + F(\cd, u, \si\pa_\o\f)
                 - F(\cd, u, Z)
           \Big)(s, B_\cd) ds 
 \\
 &=&
\int_0^\ch \big(\si\pa_\o \f -  Z\big)(s, B_\cd)\cdot \si^{-1}(s, B_\cd)dB_s
- \int_0^\ch\Big[ \frac{c}{2} 
                 + (\si\pa_\o\f -  Z)\cdot\a
           \Big](s, B_\cd) ds  \\
 &=&   
  \int_0^\ch\big(\si\pa_\o \f -  Z\big)(s, B_\cd)
             \cdot\big(\si^{-1}(s, B_\cd)dB_s -  \a_sds\big) -\frac{c}{2}\ch,
 \q \dbP\mbox{-a.s.}
\eeaa
where $|\a|\le L$. Notice that $\si^{-1} dB_t$ is a $\dbP$-Brownian motion. Applying Girsanov Theorem one sees immediately that there exists $\tilde \dbP\in \cP_L$ equivalent to $\dbP$  such that $\si^{-1}dB_t -  \a_tdt$ is a $\tilde \dbP$-Brownian motion. Then the above inequality holds $\tilde \dbP$-a.s., and by the definition of $\underline\cA^L u$:
 \beaa
 0 
 &\le& 
 \dbE^{\tilde \dbP}\big[(\f - u)_\ch\big] 
 \;\le\; -\frac{c}{2} \dbE^{\tilde \dbP}[\ch] 
 \;<\; 0,
\eeaa
which is the required contradiction.
\\
{\it Step 2.} We now show the dynamic programming principle together with the following regularity of $u$:  there exists a modulus of continuity function $\overline\rho_0$ such that,
\bea
\label{semireg}
|u(t,\o)|\le C ~\mbox{and}~ |u(t,\o) - u(t', \o')| \le C \overline\rho_0\Big(\dbf_\infty\big((t,\o), (t',\o)\big)\Big),~~ t\le t', \o,\o'\in \O.
\eea
 Indeed, by standard arguments it is clear that, for any $p\ge 1$,
\beaa
&&\dbE^{\dbP_0^t} \Big[ \|\cX^{t,\o}\|_T^p +  \|\cY^{t,\o}\|_T^p + \big(\int_t^T  |\cZ^{t,\o}_s|^2ds\big)^{p/2}\Big] \le C_p;\\
&&\dbE^{\dbP^{t,\o}} \Big[\|B^t\|_T^p +  \|Y^{t,\o}\|_T^p + \big(\int_t^T  |[\si^{t,\o}(s, B^t)]^{-1}Z^{t,\o}_s|^2ds\big)^{p/2}\Big] \le C_p;\\
&&\dbE^{\dbP_0^t} \Big[\|\cX^{t,\o} - \cX^{t,\o'}\|_T^2 \Big] \le C \rho_0(\|\o-\o'\|_t)^2;
\eeaa
and, since $\rho_0$ has polynomial growth,
\beaa
&&\dbE^{\dbP_0^t} \Big[ \|\cY^{t,\o} - \cY^{t,\o'}\|_T^2+ \int_t^T  |\cZ^{t,\o}_s - \cZ^{t,\o'}_s|^2ds\Big]\\
&\le&C \rho_0(\|\o-\o'\|_t)^2 + C\dbE^{\dbP_0^t} \Big[C \rho_0(\|\cX^{t,\o} - \cX^{t,\o'}\|_T)^2\Big] \le C\rho_1(\|\o-\o'\|_t)^2,
\eeaa
for some modulus of continuity function $\rho_1$.
In particular, this implies that
\bea
\label{semireg1}
|u(t,\o)|\le C &\mbox{and}& |u(t,\o) - u(t, \o')| \le C \rho_1(\|\o-\o'\|_t).
\eea
Given the above regularity, by standard arguments in BSDE theory, we have the following dynamic programming principle:  for any $t<t'\le T$,
 \bea\label{semiDPP}
 Y^{t,\o}_s 
\!= \!
 u^{t,\o}(t', B^t) 
 +\!\! 
 \int_s^{t'} F^{t,\o}(r, B^t_\cd,Y^{t,\o}_r, Z^{t,\o}_r)dr 
 - \!\!
 \int_s^{t'} Z^{t,\o}_r\cdot (\si^{t,\o}(r, B^t_\cd))^{-1}dB^t_r,
 \dbP^{t,\o}\mbox{-a.s.}
 \eea
In particular,  $Y^{t,\o}_s = u^{t,\o}(s, B^t)$ for all  $t\leq s\leq T$, $\dbP^{t,\o}$-a.s. That is, $Y^{t,\o}_s = u(s, \o\otimes_t B^t) = Y^{s, \o\otimes_t B^t}_s$, $\dbP^{t,\o}$-a.s.

Denote $\d := \dbf_\infty\big((t,\o), (t',\o)\big)$. 
Then
 \bea
 |u_t - u_{t'}|(\o)
 &=& 
 \Big| \dbE^{\dbP^{t,\o}}
       \Big[Y^{t,\o}_t - Y^{t,\o}_{t'} 
            + u^{t,\o}(t', B^t) - u(t', \o)
       \Big] 
 \Big|
 \nonumber\\
 &=& 
 \Big|\dbE^{\dbP^{t,\o}}
      \Big[\int_t^{t'} 
           F^{t,\o}(r,B^t_\cd,Y^{t,\o}_r, Z^{t,\o}_r)dr
           + u^{t,\o}(t', B^t) - u(t', \o)]
      \Big]
 \Big| 
 \nonumber\\
 &\le& 
 \dbE^{\dbP^{t,\o}}
 \Big[\int_t^{t'} 
      |F^{t,\o}(r,B^t_\cd,Y^{t,\o}_r, Z^{t,\o}_r)|dr 
      + C\rho_1\Big(\d+\|B^t\|_{t'}
               \Big)
 \Big],
 \label{uregularint}
\eea
Notice that
 \beaa
&& \dbE^{\dbP^{t,\o}}
 \Big[\int_t^{t'}
      \big|F^{t,\o}(r, B^t_\cd,Y^{t,\o}_r,Z^{t,\o}_r)\big|dr
 \Big] 
 \;\le\; 
 C \dbE^{\dbP^{t,\o}}
   \Big[\int_t^{t'} 
         \big(1+|Y^{t,\o}_r| + |Z^{t,\o}_r|\big)dr 
   \Big]
 \\
 &\le& 
 C\sqrt{\d} 
 \Big( \dbE^{\dbP^{t,\o}}
       \Big[\int_t^{t'}
            \big(1+|Y^{t,\o}_r|^2 + |Z^{t,\o}_r)|^2\big)dr 
       \Big]
 \Big)^{1/2} 
 \;\le\; 
 C\sqrt{\d}.
 \eeaa 
As for the second term, since $\rho_0$ has polynomial growth, one can easily see that we may assume without loss of generality that $\rho_1$ also has polynomial growth.  Note that $t'-t \le \d$. Then it is clear that there exists a modulus of continuity function $\overline\rho_0$ such that 
\beaa
 \dbE^{\dbP^{t,\o}}
 \Big[\rho_1\Big(\d + \|B^t\|_{t'}\Big)
 \Big] &\le& \overline\rho_0(\d).
 \eeaa
 Without loss of generality we assume $\overline \rho_0(\d) \ge \sqrt{\d}$. Then, plugging the last estimates into \reff{uregularint} and combining with \reff{semireg1}, we obtain \reff{semireg}.
 
 Moreover, given the regularity in $t$, we may extend the dynamic programming principle \reff{semiDPP} to stopping times, proving \reff{semiDPP1}.
\qed

\begin{rem}
\label{rem-semi1}
{\rm  For FBSDE \reff{BSDEstrong} with $(t,\o)=(0,{\bf 0})$, we have  $\cY_s := u(s, \cX_\cd)$ $\dbP_0$-a.s.  This extends the nonlinear Feynman-Kac formula of \cite{PP2}  to the path-dependent case.
\qed}
\end{rem}

\subsection{Path dependent HJB equations and 2BSDEs}

Let $\dbK$ be a measurable space (equipped with some $\si$-algebra).  We  now consider the following path dependent HJB equation:
\bea
\label{HJB}
&-\pa_t u - G(t,\o, u, \pa_\o u, \pa_{\o\o}^2 u)  =0,\q u(T,\o) = \xi(\o);&\\
&\mbox{where}~ G(t,\o, y,z,\g) := \sup_{k\in \dbK}\Big[{1\over 2} \si^2(t,\o, k): \g + F(t,\o, y, \si(t,\o, k)z, k)\Big],&\nonumber
\eea
where $\si \in \dbS^d$ and $F$ are $\dbF$-progressively measurable in all variables, and $\xi$ is $\cF_T$-measurable.  We shall assume

\begin{assum}
\label{assum-HJB}
{\rm (i)} $\si$, $F(t,\o, 0, {\bf 0},k )$, and $\xi$ are bounded by $C_0$, and $\si>{\bf 0}$.
\\
{\rm (ii)}  $\si$  is uniformly Lipschitz continuous in $\o$, and $F$ is uniformly Lipschitz contin. in $(y,z)$.
\\
{\rm (iii)} $F$ and $\xi$ are uniformly continuous in $\o$, and the common modulus of continuity function $\rho_0$ has polynomial growth. 
\\
{\rm (iv)} $\si(\cd, k)$,  $F(\cd, y,z,k)$, and $G(\cd, y,z)$ are right continuous in $(t,\o)$ under $\dbf_\infty$ for any $(y,z,k)$, in the sense of Definition \ref{defn-rightcont}.
\end{assum}

For each $t$, let $\cK^t$ denote the set of $\dbF^t$-progressively measurable $\dbK$-valued processes on $\L^t$. For any $(t,\o)\in\L$ and $k\in \cK^t$,  let $\cX^{t,\o, k}$ denote the solution to   the following SDE:
 \beaa
  \cX_s 
 = 
 \int_t^s \si^{t,\o}(r, \cX_\cd, k_r)d B^t_r,\q t\le s\le T,&&
 \dbP_0^t\mbox{-a.s.} 
 \eeaa
 Denote  $\dbP^{t,\o,k} := \dbP^t_0 \circ (\cX^{t,\o, k})^{-1}$. Since $\si>{\bf 0}$, as discussed in \cite{STZ-duality}  $\cX^{t,\o, k}$ and $B^t$ induce the same $\dbP_0^t$-augmented filtration, and thus there exists $\tilde k\in \cK^t$ such that $\tilde k (\cX^{t,\o, k}_\cd) = k$, $\dbP_0^t$-a.s.    Let  $(Y^{t,\o, k}, Z^{t,\o, k})$ denote the solution to the following BSDE on $[t, T]$:
   \beaa
  Y_s 
  = 
  \xi^{t,\o}(B^t) 
  + \int_s^T F^{t,\o}(r, B^t_\cd,  Y_r, Z_r, \tilde k_r) dr 
  - \int_s^T Z_r\cdot(\si^{t,\o}(r, B^t_\cd, \tilde k_r))^{-1}dB^t_r,
  ~\dbP^{t,\o,k}\mbox{-a.s.}
 \eeaa
We now consider the stochastic control problem:
 \beaa
 u(t,\o) := \sup_{k\in \cK^t} Y^{t,\o,k}_t,
 &(t,\o)\in\L.&
 \eeaa
We observe that this process $u$ was considered by Nutz \cite{Nutz-stochastic control}, in the stochastic control context, and shown to be the solution of a second order BSDE. The next result shows that our notion of viscosity solution is also suitable for this stochastic control problem.

\begin{prop}
\label{prop-HJB}
Under Assumption \ref{assum-HJB}, $u$ is a viscosity solution of PPDE \reff{HJB}.
\end{prop}

\proof By Proposition \ref{prop-change}, without loss of generality we assume
 \bea
 \label{HJBmonotone}
 \mbox{$G$, hence $F$, is increasing in $y$.}
 \eea
Following similar arguments as in Proposition \ref{prop-semi}, we may prove that 
\beaa
|u(t,\o)|\le C,\q  |u(t,\o) - u(t', \o')| \le C \overline\rho_0\Big(\dbf_\infty\big((t,\o), (t',\o)\big)\Big), &\mbox{for any}&(t,\o), (t', \o')\in\L.
\eeaa
This regularity, together with the standard arguments, see e.g. \cite{STZ-duality} or \cite{Peng-g}, implies further  the following dynamic programming principle:
\bea
\label{HJBDPP}
 u(t,\o) = \sup_{k\in \cK^t} \cY^{t,\o,k}_t(\t, u^{t,\o}(\t, \cd)), &\mbox{for any}&(t,\o)\in\L, \t\in \cT^t,
 \eea
where, for any $\cF^t_\t$-measurable random variable $\eta$, $(\cY, \cZ) := (\cY^{t,\o,k}(\t, \eta), \cZ^{t,\o,k}(\t, \eta))$ solves the following BSDE on $[t,\t]$:
 \beaa
 \cY_s 
 = 
 \eta(B^t_\cd) 
 + \int_s^\t F^{t,\o}(r, B^t,  \cY_r, \cZ_r, \tilde k_r) dr 
 - \int_s^\t \cZ_r\cdot(\si^{t,\o}(r, B^t_\cd, \tilde k_r))^{-1} dB^t_r,
 &&
 \dbP^{t,\o,k}-\mbox{a.s.}
 \eeaa
 We now prove the viscosity property, for the same $L$ as in Proposition \ref{prop-semi}. Again we shall only prove it  at $(t,\o) = (0, {\bf 0})$ and we will omit the superscript $^{0,{\bf 0}}$.  However, since in this case $u$ is defined through a supremum, we need to prove the viscosity subsolution property and supersolution property differently. 

\no {\it Viscosity $L-$subsolution property.} Assume to the contrary that, 
 \beaa
 c:= -\big\{\pa_t \f 
            +G(\cd, u, \pa_\o \f, \pa_{\o\o}^2 \f)
      \big\}(0,{\bf 0}) 
 \;>\; 0 
 &\mbox{for some}&
 \f\in \underline \cA^{L} u(0,{\bf 0}).
 \eeaa
As in Proposition \ref{prop-semi}, let $\ch = \ch_\e \in \cH$ be the hitting time corresponding to $\f$ in  \reff{cA}. Since $\f \in C^{1,2}(\L)$, $u\in UC_b(\L)$,  and by Assumption \ref{assum-HJB} (iv) $G$ is right continuous in $(t,\o)$  under $\dbf_\infty$, we may assume $\e$ is small enough such that 
\beaa
 -\big\{\pa_t \f 
        + G(\cd, u, \pa_\o \f, \pa_{\o\o}^2 \f)
  \big\}(t,\o) 
 \;\ge\; 
 {c\over 2} 
 \;>\;
 0,\q t\in [0, \ch].
 \eeaa
By the definition of $G$, this implies that, for any $t\in [0, \ch]$ and $k\in\dbK$,
 \beaa
  -\big\{\pa_t \f + {1\over 2} \si^2(t,\o, k):\pa_{\o\o}^2 \f 
         + F(t,\o, u, \si(\cd, k)\pa_\o \f, k)
   \big\}(t,\o) 
  \;\ge\; {c\over 2} \;>\; 0.
 \eeaa
Now for any $k \in \cK$, notice that $d \la B \ra_t = \si^2(t, B_\cd, \tilde k_t)dt$, $\dbP^k$-a.s. Denote $(\cY^k, \cZ^k) := (\cY^k(\ch, u(\ch,\cd)),  \cZ^k(\ch, u(\ch,\cd)))$. One can easily see that $u(s, B) \ge \cY^k_s$, $0\le s\le \ch$, $\dbP^k$-a.s.  For any $\d>0$, applying functional It\^{o}'s formula on $\f$ we see that, : 
\beaa
 (\f - \cY^k)_0 -(\f - u)_{\ch\wedge\d} 
 &\ge& 
 (\f - \cY^k)_0 -(\f - \cY^k)_{\ch\wedge\d} 
 \\
 &=& 
 -\int_0^{\ch\wedge\d}\Big[\pa_t \f + {1\over 2} \si^2 :\pa_{\o\o}^2\f
                   +F(\cd,\cY^k,\cZ^k)
              \Big](s, B_\cd, \tilde k_s) ds \\
              &&
 - \int_0^{\ch\wedge\d}\big(\pa_\o \f - \si^{-1}\cZ^k\big)
                           (s, B_\cd, \tilde k_s)\cdot dB_s
 \\
 &\ge&  
 \int_0^{\ch\wedge\d}\Big[ {c\over 2} 
                   + F(\cd, u, \si\pa_\o \f)
                   - F(\cd,\cY^k,\cZ^k)
             \Big](s,B_\cd,\tilde k_s) ds\\&& 
 - \int_0^{\ch\wedge\d}\big(\pa_\o \f - \si^{-1}\cZ^k\big)
               (s, B_\cd, \tilde k_s) \cdot dB_s,\q \dbP^k\mbox{-a.s.}
 \eeaa
Note again that $\cY^k_s \le u(s, B_\cd)$. Then by \reff{HJBmonotone} we have
 \beaa
 &&\big(u - \cY^k\big)_0 -\big(\f - u\big)_{\ch\wedge\d} = \big(\f - \cY^k\big)_0 -\big(\f - u\big)_{\ch\wedge\d}
 \\
 &\ge&  
 \int_0^{\ch\wedge\d}\Big[{c\over 2} 
                  + F(\cd, u, \si\pa_\o \f)
                  - F(\cd, u, \cZ^k)
             \Big](s, B_\cd, \tilde k_s) ds \\
&& - \int_0^{\ch\wedge\d}\big(\pa_\o \f - \si^{-1}\cZ^k\big)
               (s, B_\cd, \tilde k_s)\cdot dB_s
 \\
 &=&
 \int_0^{\ch\wedge\d}\Big[{c\over 2} 
                  + (\si\pa_\o \f -  \cZ^k)\cd \a
             \Big](s, B_\cd,\tilde k_s) ds 
 - \int_0^{\ch\wedge\d}\big(\pa_\o \f - \si^{-1}\cZ^k\big)
               (s, B_\cd,\tilde k_s) \cdot dB_s
 \\
 &=&
 {c\over 2}(\ch\wedge \d)
 - \int_0^{\ch\wedge\d}\big(\si\pa_\o \f - \cZ^k\big)
               (s, B_\cd,\tilde k_s)\cdot (\si^{-1}(s, B_\cd, \tilde k_s)dB_s - \a_sds),
 \q \dbP^k\mbox{-a.s.}
\eeaa
where $|\a|\le L$ and $\l$ is bounded.  As in Proposition \ref{prop-semi}, we may define $\tilde \dbP^{k}\in \cP_L$ equivalent to $\dbP$  such that $\si^{-1}(t, B_\cd, \tilde k_t)dB_t - \a_tdt$ is a $\tilde \dbP^{k}$-Brownian motion. Then the above inequality holds $\tilde\dbP^{k}$-a.s., and by the definition of $\underline\cA^L u$, we have
 \beaa
 u_0 - \cY^{k}_0  
 \ge
 u_0 - \cY^{k}_0 
 - \dbE^{\tilde \dbP^{k}}\big[(\f - u)_{\ch\wedge \d}\big] \ge
 {c\over 2} \dbE^{\tilde \dbP^{k}}[\ch\wedge \d] \ge {c\over 2} \d\Big[1- \tilde \dbP^{k}[ \ch \le \d]\Big].
 \eeaa
By \reff{chest}, for $\d$ small enough we have
 \beaa
 u_0 - \cY^{k}_0  
 &\ge& {c\over 2} \d \Big[1-  C\e^{-4} \d^2\Big] \ge {c\d\over 4}>0.
 \eeaa
This implies that $u_0 - \sup_{k\in\cK} \cY^k_0 \ge {c\d\over 4} >0$, which is in contradiction with \reff{HJBDPP}.
\\
{\it Viscosity $L-$supersolution property.} Assume to the contrary that, 
 \beaa
 c
 := 
 \Big\{\pa_t \f 
       + G(\cd, u, \pa_\o \f, \pa_{\o\o}^2 \f)\Big](0,{\bf 0}) 
 \;>\; 0 
 &\mbox{for some}&\f\in \overline \cA^{L} u(0,{\bf 0}).
\eeaa
By the definition of $F$, there exists $k_0\in \dbK$ such that
 \beaa
 \Big\{\pa_t \f + {1\over 2} \si^2(\cd, k_0):\pa_{\o\o}^2 \f 
       +F(\cd, u, \si(\cd, k_0)\pa_\o \f, k_0)
 \Big\}(0, {\bf 0}) 
 \;\ge\; {c\over 2} 
 \;>\; 0
\eeaa
Again, let $\ch = \ch_\e \in \cH$ be the hitting time corresponding to $\f$ in  \reff{cA}, and by the right continuity of $\si$ and $F$ in Assumption \ref{assum-HJB} (iv) we may assume $\e$ is small enough so that 
 \beaa
  \Big\{\pa_t \f + {1\over 2} \si^2(\cd, k_0):\pa_{\o\o}^2 \f 
        +F(\cd, u, \si(\cd, k_0)\pa_\o \f, k_0)
  \Big\}(t,\o) 
 \;\ge\; 
 {c\over 3} \;>\;0,\q t\in [0, \ch].
 \eeaa
 Consider the constant process $k:= k_0 \in \cK$. It is clear that the corresponding $\tilde k = k_0$. Follow similar arguments as in the subsolution property, we arrive at the following contradiction:
$
u_0 - \cY^{k}_0  \le- {c\over 3} \dbE^{\tilde \dbP^{k}}[\ch] < 0.
$
\qed

\begin{eg}
\label{eg-2BSDE}
{\rm Assume $\dbK := \{k \in \dbS^d: \underline \si \le k \le \overline\si\}$, where ${\bf 0}<\underline \si <\overline\si$ are constant matrices. Set $\si(t,\o, k) := k$. Then $Y_t(\o) = u(t,\o)$ is the solution to the following second order BSDE, as introduced by  \cite{STZ-2BSDE}:
 \bea
 \label{2BSDE}
 Y_t 
 = 
 \xi(B_\cd) 
 + \int_t^T F(s, B_\cd, Y_s, Z_s, \hat a_s^{1\over 2})ds 
 - \int_t^T Z_s \cdot (\hat a_s)^{-{1\over 2}} dB_s 
 - d K_t, 
 \cP\mbox{-q.s.}
 \eea
where $\cP := \{ \dbP\in \cP^0_\infty: \a^\dbP=0, \b^\dbP \in \dbK\}$, $\hat a$ is the universal process such that $d \la B\ra_t = \hat a_t dt$, $\cP$-q.s. and $K$ is an increasing process satisfying certain minimum condition. 
\qed}
\end{eg}

\begin{rem}
\label{rem-Isaacs}
{\rm By using the zero-sum game, we may also obtain a representation formula for the viscosity solution of the following path dependent Bellman-Isaacs equation:
\bea
\label{HJBI}
-\pa_t u - G(t,\o, u, \pa_\o u, \pa_{\o\o}^2 u)  =0,\q u(T,\o) = \xi(\o), 
\eea
where
\beaa
G(t,\o, y,z,\g) 
&:=& 
\sup_{k_1\in \dbK_1}\inf_{ k_2\in \dbK_2}
\Big[{1\over 2} \si^2(t, k_1, k_2): \g 
     + F(t,\o, y, \si(t, k_1, k_2)z, k_1, k_2)\Big]
\\
&=&
\inf_{ k_2\in \dbK_2}\sup_{k_1\in \dbK_1}
\Big[{1\over 2} \si^2(t,k_1, k_2): \g 
     + F(t,\o, y, \si(t, k_1, k_2)z, k_1, k_2)\Big].
\eeaa
See Pham and Zhang \cite{PZ}.
\qed}
\end{rem}

\section{Stability and Partial Comparison}
\label{sect-stability}
\setcounter{equation}{0}

\subsection{Stability}
The main result of this section is the following extension of Theorem 4.1 in \cite{EKTZ}, with a proof following the same line of argument. However the present fully nonlinear context makes a crucial use of Theorem \ref{thm-optimal}.  Denote, for any $(t,y, z, \g) \in [0, T) \times \dbR\times \dbR^d\times \dbS^d$ and $\d>0$,
\beaa
O_\d(t, y, z, \g) &:=& \Big\{ (s, \tilde \o, \tilde y, \tilde z, \tilde \g)\in \L^t \times \dbR\times \dbR^d\times \dbS^d: \nonumber\\
&&\dbf^t_\infty((s,\tilde\o), (t, {\bf 0})) + |\tilde y-y|+|\tilde z-z| + |\tilde \g-\g|\le \d\Big\}.
\eeaa

\begin{thm}
\label{thm-stability}  Let $L>0$, $G$ satisfy Assumption \ref{assum-G}, and $u\in \Usub$ (resp. $u\in \Usup$). Assume

(i) for any $\e>0$, there exist $G^\e$ and $u^\e\in \Usub$ (resp. $u^\e\in \Usup$) such that $G^\e$ satisfies Assumption \ref{assum-G} and $u^\e$ is a viscosity $L$-subsolution (resp. $L$-supersolution) of PPDE \reff{PPDE} with generator $G^\e$;

(ii) as $\e\to 0$, $(G^\e, u^\e)$ converge to $(G, u)$  locally uniformly in the following sense: for any $(t,\o, y, z, \g)\in \L\times \dbR\times \dbR^d\times \dbS^d$, there exists $\d>0$ such that, 
\bea
\label{localuniform}
\lim_{\e\to 0} \sup_{(s, \tilde \o, \tilde y, \tilde z, \tilde \g) \in O_\d(t, y, z, \g)}\Big[|(G^\e - G)^{t,\o}(s,  \tilde \o, \tilde y, \tilde z, \tilde \g) |+ |(u^\e - u)^{t,\o} (s,  \tilde \o)| \Big] = 0.
\eea
Then $u$ is a viscosity $L$-subsolution (resp. $L$-supersol.) of PPDE \reff{PPDE} with generator $G$.
\end{thm}
\proof Without loss of generality we shall only prove the viscosity subsolution property at $(0, {\bf 0})$. Let $\f\in {\overline\cA}^{L}\!u(0,{\bf 0})$  with corresponding $\ch \in \cH$, $\d_0>0$ be a constant  such that $\ch_{\d_0} \le \ch$ and $\lim_{\e\to 0} \rho(\e,\d_0)=0$, where
\beaa
&\dis \rho(\e,\d ) :=  \sup_{(t, \o, \tilde y, \tilde z, \tilde \g) \in O_{\d}(0, y_0, z_0, \g_0) }\Big[|G^\e - G|(t, \o, \tilde y, \tilde z, \tilde \g) + |u^\e - u| (t, \o)\Big],&\\
&\mbox{and}~~( y_0, z_0, \g_0):=( \f_0, \pa_\o \f_0, \pa_{\o\o}^2\f_0).&
\eeaa 
For $0 <\d \le \d_0$, denote $\f_\d(t,\o) := \f(t,\o) + \d t$. By \reff{cA} and Lemma \ref{lem-che} we have
\beaa
(\f_\d-u)_0= (\f-u)_0 = 0 \le \underline \cE^L_0\Big[(\f-u)_{\ch_\d}\Big] <  \underline \cE^L_0\Big[(\f_\d -u)_{\ch_\d}\Big].
\eeaa 
By \reff{localuniform}, there exists $\e_\d>0$ small enough such that, for any $\e \le \e_\d$, 
\bea
\label{stabilityest1}
(\f_\d-u^\e)_0  <  \underline \cE^L_0\Big[(\f_\d -u^\e)_{\ch_\d}\Big].
\eea
Denote $X := X^{\e,\d} := u^\e - \f_\d \in \Usub$.  Define $\hat X := X \1_{[0, \ch_\d)} + X_{\ch_\d-} \1_{[\cd_\d, T]}$, $Y := \overline \cE^L[\hat X]$, and $\t^* := \inf\{t \ge 0: Y_t = \hat X_t\}$,  as in Theorem \ref{thm-optimal}.  Then all the results in   Theorem \ref{thm-optimal} hold.  Noticing that $X_{\ch_\d-} \le X_{\ch_\d}$, by \reff{stabilityest1}  we have
\beaa
\overline\cE^L_0[\hat X_{\ch_\d}] \le \overline\cE^L_0[X_{\ch_\d}] = - \underline \cE^L_0\Big[(\f_\d -u^\e)_{\ch_\d}\Big] < - (\f_\d-u^\e)_0 = X_0 \le Y_0 =  \overline\cE^L_0[Y_{\t^*}] = \overline\cE^L_0[\hat X_{\t^*}]. 
\eeaa
Then there exists $\o^*$ such that $t^*:= \t^*(\o^*) < \ch_{\d}(\o^*)$, and thus $\ch_\d^{t^*, \o^*} \in \cH^{t^*}$. We shall remark though that here $Y, \t^*, \o^*, t^*$ all depend on $\e, \d$. Now define 
\beaa
\f^\e_\d(t, \o) := \f_\d^{t^*, \o^*}(t,\o) - \f_\d(t^*, \o^*) + u^\e(t^*, \o^*),\q (t,\o) \in \L^{t^*}. 
\eeaa
It is straightforward to check that $\f^\e_\d \in {\underline\cA}^{L} u^\e(t^*, \o^*)$ with corresponding hitting time $\ch_\d^{t^*, \o^*}$. Since $u^\e$ is a viscosity $L$-subsolution of PPDE \reff{PPDE} with generator $G^\e$, we have 
\bea
\label{stabilityest2}
0 &\ge& \Big[-\pa_t \f^\e_\d - (G^\e)^{t^*,\o^*}(\cd, \f^\e_\d, \pa_\o \f^\e_\d, \pa_{\o\o}^2 \f^\e_\d)\Big](t^*, {\bf 0})\nonumber\\
&=&  \Big[-\pa_t \f - \d - G^\e(\cd, u^\e, \pa_\o \f, \pa_{\o\o}^2 \f)\Big](t^*, \o^*).
\eea

Note that $t^* < \ch_\d(\o^*)$, then $|u^\e- u|(t^*, \o^*)  \le \rho(\e, \d)\le \rho(\e, \d_0)$. By \reff{localuniform} and Definition \ref{defn-rightcont}, we may set $\d$ small enough and then $\e$ small enough so that $(\cd, u^\e, \pa_\o \f, \pa_{\o\o}^2 \f)(t^*, \o^*) \in O_{\d_0}(0, y_0, z_0, \g_0)$. Thus, \reff{stabilityest2} leads to
\beaa
0 &\ge &  \Big[-\pa_t \f - \d - G^\e(\cd, u^\e, \pa_\o \f, \pa_{\o\o}^2 \f)\Big](t^*, \o^*) \\
&\ge&  \Big[-\pa_t \f  - G(\cd, u^\e, \pa_\o \f, \pa_{\o\o}^2 \f)\Big](t^*, \o^*) - \d - \rho(\e, \d_0)\\
&\ge&  \Big[-\pa_t \f  - G(\cd, u, \pa_\o \f, \pa_{\o\o}^2 \f)\Big](t^*, \o^*) - \d - \rho(\e, \d_0) - C \rho(\e, \d)\\
&\ge& \cL\f_0 - C\sup_{(t, \o): t < \ch_\d(\o)} \Big[  |u(t,\o) - u_0|+|\pa_\o\f(t,\o) - \pa_\o\f_0| + |\pa^2_{\o\o}\f(t,\o)-\pa^2_{\o\o}\f_0|\Big]  \\
&&-\sup_{(t, \o): t < \ch_\d(\o)} \Big|G(t,\o,y_0, z_0, \g_0) - G(0, {\bf 0},y_0, z_0, \g_0)\Big|-  \d - \rho(\e, \d_0) - C  \rho(\e, \d),
\eeaa
where we used the fact that $G$ satisfies Assumption \ref{assum-G}.  Notice that the right-continuity of $G$ in $(t,\o)$ under $\dbf_\infty$ allows us to control the last line. Now by first sending $\e\to 0$ and then $\d\to 0$ we obtain $\cL\f_0  \le 0$. Since  $\f\in {\overline\cA}^{L}\!u(0,{\bf 0})$ is arbitrary, we see that $u$ is a viscosity subsolution of PPDE \reff{PPDE} with generator $G$ at $(0, {\bf 0})$ and thus complete the proof.
 \qed
 
 \begin{rem}
\label{rem-stability}
{\rm Similar to Theorem 4.1 in \cite{EKTZ}, we need the same $L$ in the proof of Theorem \ref{thm-stability}. If $u^\e$ is only a viscosity subsolution of PPDE \reff{PPDE} with generator $G^\e$, but with possibly different $L_\e$, we are not able to show that $u$ is a viscosity subsolution of PPDE \reff{PPDE} with generator $G$.
\qed}
\end{rem}

\subsection{Partial comparison of viscosity solutions}

In this section, we prove a partial comparison principle, i.e. a comparison result of a viscosity super- (resp. sub-) solution and a classical sub- (resp. super-) solution. The proof is also crucially based on Theorem \ref{thm-optimal}. Moreover, this result is a first key step for our comparison principle in the accompanying paper \cite{ETZ2}.

\begin{prop}
\label{prop-comparison}
Let Assumption \ref{assum-G} hold true. Let $ u^1\in \Usub$ be a viscosity subsolution and $ u^2\in\Usup$ a viscosity supersolution of PPDE (\ref{PPDE}). If $ u^1(T,\cd) \le u^2(T,\cd)$ and either $ u^1$ or $ u^2$ is in $C^{1,2}(\L)$, then  $ u^1\le u^2$ on $\L$.
\end{prop}

\proof We shall only prove $u^1_0 \le u^2_0$. The inequality for general $t$ can be proved similarly. Without loss of generality, we assume $ u^1$  is a viscosity $L$-subsolution and $ u^2\in C^{1,2}(\L) $ is  a classical  $L$-supersolution. 
By Proposition \ref{prop-change}, we may assume that 
 \bea\label{Gmonotone}
 \mbox{$G$ is nonincreasing in $y$}. 
 \eea
 
Assume to the contrary that $c := {1\over 2T}[ u^1_0 - u^2_0] >0$. 
Denote 
\beaa
X_t := ( u^1 -  u^2)^+_t+ c t,\q \wh X_t := X_t \1_{\{t<T\}} + X_{T-}\1_{\{t=T\}};\q Y_t(\o) :=  \sup_{\t\in \cT^t}
 \overline\cE^L_t[\wh X^{t,\o}_\t].
 \eeaa
Since $u^1\in\Usub$ is bounded from above and $u^2\in  \Usup$ is bounded from below, it follows that $X$ is a bounded process in $\Usub$. Let $\t^*:=\inf\{t \ge 0:Y_t =\wh X_t\}$. Note that $u^1_T\le u^2_T$ and  $X\in\Usub$ has positive jumps. Then it follows from Theorem \ref{thm-optimal} that
 \beaa
 \overline\cE^L_0[\wh X_T] \le  \overline\cE^L_0[X_T] = cT < 2cT = X_0 = \wh X_0 \le Y_0 =  \overline\cE^L_0[Y_{\t^*}] =  \overline\cE^L_0[\wh X_{\t^*}].
 \eeaa
This implies that $t^* := \t^*(\o^*) <T$ for some $\o^*\in\O$. Note that 
 \beaa
 ( u^1- u^2)^+(t^*,\o^*) + c t^* 
 \;=\; 
\wh X_{t^*}(\o^*) 
 \;=\; 
 Y_{t^*}(\o^*) 
 \;\ge\; 
 \overline\cE^L_{t^*}[X_{T-}^{t^*,\o^*}] 
 \;\ge\; 
 cT
 \;>\;
 0.
 \eeaa
Then $( u^1- u^2)(t^*, \o^*) >0 $. Since $ (u^1- u^2)^{t^*,\o^*}\in \Usub^{t^*}$, there exists $\ch \in \cH^{t^*}$ such that $\ch < T$ and $( u^1- u^2)^{t^*, \o^*}_t >0$ for all $t\in [t^*, \ch]$, and thus $\wh X^{t^*, \o^*}_t = X^{t^*, \o^*}_t = ( u^1 -  u^2)^{t^*, \o^*}_t + c t$ for all $t\in [t^*, \ch]$.

Now observe that $\f(t,\o):= (u^2)^{t^*, \o^*}(t,\o) - c t + X_{t^*}(\o^*)$ $\in$ $C^{1,2}(\L^{t^*})$, a consequence of our assumption $ u^2\in C^{1,2}(\L)$. Moreover, for any $\t \in \cT^{t^*}$, it follows from the $\overline\cE^L$-supermartingale property of the nonlinear Snell envelope $Y$ that
 \beaa
 \big( (u^1)^{t^*, \o^*}-\f\big)_{t^*}  &=& 0 \;=\; Y_{t^*}(\o^*) - X_{t^*}(\o^*)
 \;\ge\;
 \overline\cE^L_{t^*}\big[Y^{t^*,\o^*}_{\t\wedge\ch}\big]  - X_{t^*}(\o^*)\\
 &\ge&
 \overline\cE^L_{t^*}\big[X^{t^*,\o^*}_{\t\wedge\ch} \big] - X_{t^*}(\o^*)
 \;=\;
 \overline\cE^L_{t^*}
 \big[\big( (u^1)^{t^*, \o^*}-\f\big)_{\t\wedge\ch}\big].
 \eeaa
By the arbitrariness of $\t\in\cT^{t^*}$, and the fact that $\underline{\cE}^L[\cdot]=-\overline{\cE}^L[-\cdot]$, this proves that $\f \in \underline\cA^L u^1(t^*, \o^*)$, and by the viscosity $L$-subsolution property of $ u^1$:
 \beaa
 0 
 &\ge&
 \big\{-\pa_t \f 
       - G(.,  u^1, \pa_\o\f, \pa^2_{\o\o}\f)\big\}
 (t^*,\o^*)\\
 &=& 
 c-
 \big\{\pa_t u^2
       + G(.,  u^1,\pa_\o u^2,
                           \pa^2_{\o\o} u^2)\big\}
 (t^*,\o^*)
 \\
 &\ge& 
 c-
 \big\{\pa_t u^2
       + G(.,  u^2,\pa_\o u^2,
                           \pa^2_{\o\o} u^2)\big\}
 (t^*,\o^*),
 \eeaa
where the last inequality follows from \reff{Gmonotone}. Since $c>0$, this is in contradiction with the supersolution property of $ u^2$. 
 \qed

As a direct consequence of the above partial comparison, we have
\begin{prop}
\label{prop-unique}
Let Assumption \ref{assum-G} hold true. If PPDE \reff{PPDE} has a classical solution $u\in C^{1,2}(\L)\cap UC_b(\L)$, then $u$ is the unique viscosity solution of PPDE \reff{PPDE} with terminal condition $u(T,\cd)$.
\end{prop} 

In our accompanying paper \cite{ETZ2}, we shall prove the uniqueness of viscosity solutions without assuming the existence of classical solutions.

\section{Viscosity Solutions of  Backward Stochastic PDEs}
\label{sect-BSPDE}
\setcounter{equation}{0}

In this section, we show that our PPDEs includes Backward SPDEs as a special case.  We remark that such BSPDEs arise naturally in many applications, see e.g.  \cite{MYZ} and \cite{OSZ}.  Consider the following BSPDE with $\dbF$-progressively measurable solution $(u, q)$:
\bea
\label{BSPDE}
u(t,\o,x) 
= 
\xi(\o, x) + \int_t^T F(s,\o,x, u, Du, D^2u, q, Dq) ds 
- \int_t^T q(s,\o,x)\cdot dB_s,\dbP_0\mbox{-a.s.}
\eea
where $x\in\dbR^{d'}$, and $D,D^2$ denote the gradient and Hessian with respect to the $x-$variable.
Assume $ u\in C^{1,2}(\L \times \dbR^{d'})$, namely $\pa_t u, \pa_\o u, D u, \pa_{\o\o}^2u, D \pa_\o u, D^2 u$ exist and are continuous, where the derivatives in $x$ are in standard sense and the smoothness in $(t,\o)$ is in the sense of Definition \ref{defn-spaceC12}.   Fix $x$ and apply funtional  It\^{o}'s formula, we have
 \beaa
 d u (t,\o,x) 
 = 
 \big(\pa_t u + {1\over 2} \tr(\pa^2_{\o\o}u) \big)(t,\o,x)dt 
 + \pa_{\o} u(t,\o,x)\cdot dB_t,
 ~~\dbP_0\mbox{-a.s.}
 \eeaa
Comparing this with \reff{BSPDE} we obtain
 $$
 q(t,\o,x) = \pa_{\o}u(t,\o,x),
 ~\big(\pa_t u + {1\over 2}\tr(\pa^2_{\o\o}u)\big)(t,\o,x) 
    + F(t,\o,x,u, Du, D^2u, q, Dq)=0.
 $$
This leads to a mixed PPDE:
 \bea\label{BSPDE-o}
 \wh\cL u (t,\o,x)=0,
 &&  u(T,\o,x) = \xi(\o,x),\q x\in\dbR^{d'},
\eea
where, for $\f\in  C^{1,2}(\L \times \dbR^{d'})$,
 \bea\label{cL}
 \wh\cL\f
 &:=& 
 -\pa_t \f - {1\over 2}\tr(\pa^2_{\o\o}\f )
 -F(., \f, D\f, D^2\f, \pa_{\o}\f, D\pa_{\o}\f).
 \eea

To incorporate the mixed PPDE \reff{BSPDE-o} into our framework, we enlarge the space of paths to $\hat\O:=\O\times\{\o'\in C^0([0,T],\dbR^{d'}): \o'_0={\bf 0}\}$.  Denote $\hat\L:= [0,T]\times\hat\O$, and 
 \beaa
 \hat G^x(t,\hat\o, y, z,\g) 
 :=
 {1\over 2}\g_{22} + F\big(t,\o,x+\o'_t, 
                             y,z_1,\g_{11},z_2,\g_{12}\big),
 &&
 \hat\xi^x(\hat\o) 
 := 
 \xi(\o,x+\o'_T),
\eeaa
for all $x\in\dbR^{d'}$ and $(t,\hat\o, y, z,\g) \in \hat\L \times \dbR\times \dbR^{d+d'}\times \dbS^{d+d'}$. Note that  $\hat G^x(t,.)$ and $\hat\xi^x$ depend on $\hat\o=(\o,\o')$ only through the pair $(\o,\o'_t)$ and $(\o,\o'_T)$, respectively.

\begin{defn}\label{defn-BSPDE}
We say $u$ is a viscosity solution (resp. supersolution, subsolution) of BSPDE \reff{BSPDE} if, for any fixed $x$, the process $\hat u^x(t,\hat\o) := u(t,\o, x+\o'_t)$, $t\in[0,T]$, $\hat\o=(\o,\o')\in\hat\O$, is a viscosity solution (resp. supersolution, subsolution) of the PPDE:
 \beaa
 -\pa_t\hat u^x (t,\hat\o)  
 - \hat G^x\big(t,\hat\o, \hat u^x, \pa_{\hat\o}\hat u^x, 
             \pa_{\hat\o\hat\o} \hat u^x\big) 
 =0,
 &\mbox{on}~~\hat\L,~~\mbox{and}&
 \hat u^x(T,\hat\o) = \hat\xi^x(\hat\o). 
 \eeaa
\end{defn}

\begin{rem}
\label{rem-BSPDE-PDE}
{\rm When $\xi$ and $F$ do not depend on $\o$, one can easily see that $q=0$ and $u = u(t,x)$ is deterministic. Then  BSPDE \reff{BSPDE} reduces to a standard PDE. In this case, our Definition \ref{defn-BSPDE} is not the same as the standard viscosity solution of PDEs, but in the sense of Remark \ref{rem-markovian-viscosity} (i).
\qed}
\end{rem}

\begin{rem}
\label{rem-SPDE}
{\rm In the same manner we may also transform the following (forward) Stochastic PDE into a (forward) PPDE:
\bea
\label{SPDE}
u(t,\o, x)=u_0(x)+\int_0^t F(s,\o, x,u,u_x, u_{xx})ds
             +\int_0^t \si(s,\o, x,u,u_x)dB_s.
\eea
Due to its forward nature, the definition of viscosity solutions will be quite different. However, the approach which will be specified in next section and in \cite{ETZ2} still works in this case. See Buckdahn, Ma and Zhang \cite{BMZ2}.
\qed}
\end{rem}

\section{A revisit of semi-linear PPDEs}
\label{sect-semilinear}
\setcounter{equation}{0}

In \cite{EKTZ}, we proved the comparison principle for semilinear PPDE \reff{semiPPDE}, in the case $\si = I_d$.  One important argument was the Bank-Baum approximation in \cite{BB}, which unfortunately does not seem to be extendable to the fully nonlinear case. In this section we provide an alternative proof of the comparison principle for semilinear PPDE \reff{semiPPDE}.  This approach works in fully nonlinear case as well, but with much more involved technicalities, see our accompanying paper \cite{ETZ2}.  It has also been applied by Henry-Labordere, Tan,  and Touzi \cite{HTT}  to study a new type of numerical methods for BSDEs.
 
In order to focus on the main idea and simplify the presentation, we restrict to the case $\si = I_d$. That is, we shall consider the following PPDE:
\bea
\label{semiPPDE2}
\cL u(t,\o) := -\pa_t u - \frac12 I_d:\pa^2_{\o\o} u 
- F(t,\o, u, \pa_\o u)=0. 
\eea
We first give an alternative definition for viscosity solutions of semilinear PPDE \reff{semiPPDE2}. We remark that the key point in \reff{cA} and Definition \ref{defn-viscosity} is that the class $\cP_L$ covers all the probability measures induced by the linearization of the generator $G$. In the semilinear case, since the diffusion term $\si$ is already fixed, we shall only consider the drift uncertainty induced by the linearization of generator $F$, as we did in \cite{EKTZ}. To be precise, define
\beaa
 M^{t,\a}_T &:=& \exp\Big(\int_t^T \a_s\cdot dB^t_s - \frac12 \int_t^T |\a_s|^2 ds\Big);\\
\cP^{t}_L &:=& \Big\{ \dbP(\cd) := \int_\cd M^{t,\a}_T d\dbP^{t}_0: ~\a \in \dbL^0(\L^t, \dbR^d) ~\mbox{such that}~|\a|\le L\Big\};\\
\overline \cE_{t}^L[\xi] &:=& \sup_{ \dbP\in \cP^{t}_L} \dbE^{\dbP}[\xi],\q 
\underline \cE_{t}^L[\xi] \;:=\; \inf_{ \dbP\in \cP^{t}_L} \dbE^{\dbP}[\xi];
\eeaa
and
\beaa
 \underline\cA^{L}\!u(t,\o) 
 :=
 \Big\{\f\in C^{1,2}(\L^{\!t}):
       (\f-u^{t,\o})_t
       = 0=
      \inf_{\t\in \cT^t}\underline\cE_{t}^L\big[(\f-u^{t,\o})_{\cdot\wedge\ch}
                       \big]
      ~\mbox{for some}~\ch\in \cH^t\Big\},
 \\
 \overline\cA^{L}\!u(t,\o) 
 := 
 \Big\{\f \in C^{1,2}(\L^{\!t}):
      (\f-u^{t,\o})_t
      =0=
      \sup_{\t\in\cT^t}
      \overline\cE_{t}^L\big[(\f-u^{t,\o})_{\cdot\wedge\ch}
                       \big] 
      ~\mbox{for some}~\ch\in \cH^t
 \Big\}.
\eeaa

\begin{defn}
\label{defn-semi-viscosity}
\no {\rm (i)} Let $L>0$. We say $u\in\Usub$ (resp. $\Usup$) is a viscosity $L$-subsolution (resp. $L$-supersolution) of semilinear PPDE \reff{semiPPDE2}  if,  for any $(t,\o)\in [0, T)\times \O$ and any $\f \in \underline\cA^{L}u(t,\o)$ (resp. $\f \in \overline\cA^{L}u(t,\o)$):
 \beaa
 \big\{-\pa_t \f   - \frac12 I_d: \pa^2_{\o\o}\f
       - F^{t,\o}(., \f,\pa_\o \f)  
 \big\}(t,{\bf 0}) 
 &\le  ~~(\mbox{resp.} \ge)&  0.
 \eeaa

\no {\rm (ii)} $u\in\Usub$ (resp. $\Usup$)  is a viscosity subsolution (resp. supersolution) of PPDE  \reff{semiPPDE2} if  $u$ is viscosity $L$-subsolution (resp. $L$-supersolution) of PPDE \reff{semiPPDE2} for some $L>0$.

\no {\rm (iii)} $u\!\in\! UC_b(\L)$ is viscosity solution of PPDE  \reff{semiPPDE2} if it is viscosity sub- and supersolution.
\end{defn}

Under Assumption \ref{assum-F}, following almost the same arguments  and after obvious modifications when necessary, one can easily check that  Theorems \ref{thm-consistency}, \ref{thm-stability}, and Proposition \ref{prop-semi} still hold. Moreover, we may improve the partial comparison principle of Proposition \ref{prop-comparison} as follows.  First, we extend the space $C^{1,2}(\L)$:
\begin{defn}
\label{defn-barC}
Let $t\in [0,T)$, $u\in \dbL^0(\L^t)$. We say $u\in \overline C^{1,2}(\L^t)$  if  there exist  random times $t=\ch_0 \le \ch_1 \cds\le T$ (not necessarily hitting times)  such that,
\\
{\rm (i)} $\ch_i < \ch_{i+1}$ whenever $\ch_i <T$, and for all $\o\in\O^t$, the set  $\{i: \ch_i(\o) < T\}$ is finite;
\\
{\rm (ii)} For each $i$, $\o\in\O^t$, and $s\in [\ch_i(\o), \ch_{i+1}(\o))$, there exist $\ch \in \cH^s$ and $\tilde u^{s,\o} \in C^{1,2}(\L^s)$ such that $\ch < \ch^{s,\o}_{i+1}$ and $u^{s,\o} = \tilde u^{s,\o}$ on $[s, \ch]$.
\\
{\rm (iii)} 
$u$ is bounded and continuous in $t$.
\end{defn}

Roughly speaking, $\overline C^{1,2}(\L)$ consists of processes $u$ which are piecewise $C^{1,2}(\L)$, in the sense that $u$ is smooth on $[\ch_i, \ch_{i+1})$ mentioned above.

For $u\in \overline C^{1,2}(\L^t)$ and $(s,\o) \in [t, T)\times \O^t$, we may define the derivatives of $u$ at $(s,\o)$ as the derivatives of $\tilde u^{s,\o}$ at $(s, {\bf 0})$, where $\tilde u^{s,\o}$ is defined in Definition \ref{defn-barC} (ii). Clearly these derivatives are independent of the choices of $\tilde u$. We remark that the processes in  $\overline C^{1,2}(\L^t)$ are in general not continuous in $\o$. 
We also note that the space $\overline C^{1,2}(\L^t)$ here is slightly different from that in \cite{EKTZ}, and in \cite{ETZ2} we shall modify it slightly further for technical reasons.

We first extend the partial comparison principle Proposition \ref{prop-comparison} to the case that either $u^1$ or $u^2$ 
is only in $\overline C^{1,2}(\L)$, instead of $C^{1,2}(\L)$.  Our proof relies heavily on the theory of Reflected BSDEs, for which we refer to El Karoui et al \cite{EKPPQ}, Hamad\`{e}ne \cite{Hamadene}, and Peng and Xu \cite{PX}  for details. We note that Remark 3.11 in our earlier paper \cite{EKTZ} on this issue is heuristic. The precise statements are given  below. 

\begin{rem}
\label{rem-RBSDE}
{\rm Let  $X\in \dbL^0(\L)$ be  bounded and {\cad} with positive jumps. Fix $L>0$. 

(i) Let $(\tilde Y, \tilde Z, \tilde K)$ denote the unique $\dbF$-measurable solution to the following RBSDE:
\bea
\label{RBSDE}
\left\{\ba{lll}
\dis \tilde Y_t 
= X_T  + \int_t^T L|\tilde Z_s|ds 
       - \int_t^T \tilde Z_s\cdot dB_s + \tilde K_T - \tilde K_t;\\
\dis \tilde Y_t \ge X_t,\q [\tilde Y_{t-} - X_{t-}]d\tilde K_t =0;
\ea\right. \dbP_0\mbox{-a.s.}
\eea
Then $\tilde Y$ and $\tilde K$ are continuous in $t$, $\dbP_0$-a.s. Moreover, there exists $\t^*\in \cT$ such that  
\bea
\label{tau*}
\t^* = \inf\{t\ge 0: \tilde Y_t = X_t\},\q \dbP_0\mbox{-a.s.}
\eea
We remark that, since we require $\tilde Y$ to be $\dbF$-measurable, we cannot claim $\tilde Y$ is continuous for all $\o$.  Consequently, the right side of \reff{tau*} may not be an $\dbF$-stopping time, but a stopping time adapted to the $\dbP_0$-augmented filtration of $\dbF$.

 (ii) Define
\beaa
Y_t(\o) &=& \sup_{\t\in \cT^t} \overline\cE^L_t[X^{t,\o}_\t],\q (t,\o) \in \L.
\eeaa
Then $Y_0 = \tilde Y_0$.  Moreover, for any $\t\in \cT$, following standard arguments one may easily show that, for $\dbP_0$-a.e. $\o$,  $(\tilde Y^{\t, \o}, \tilde Z^{\t, \o}, \tilde K^{\t, \o})$ satisfies the following RBSDE on $[\t(\o), T]$:
 \beaa
\left\{\ba{lll}
\dis \tilde Y^{\t, \o}_t 
= 
X^{\t, \o}_T  + \int_t^T L|\tilde Z^{\t, \o}_s|ds 
              - \int_t^T \tilde Z^{\t, \o}_s\cdot dB^{\t(\o)}_s 
              + \tilde K^{\t, \o}_T - \tilde K^{\t, \o}_t;
\\
\dis \tilde Y^{\t, \o}_t \ge X^{\t, \o}_t,\q [\tilde Y^{\t, \o}_{t-} - X^{\t, \o}_{t-}]d\tilde K^{\t, \o}_t =0;
\ea\right. \dbP^{\t(\o)}_0\mbox{-a.s.}
\eeaa
Then, $Y_{\t(\o)} (\o)= \tilde Y^{\t, \o}_{\t(\o)} = \tilde Y_{\t(\o)}(\o)$. That is, $Y_\t = \tilde Y_\t$, $\dbP_0$-a.s. for all $\t \in \cT$. In other words, $Y$ and $\tilde Y$ are $\dbP_0$-modifications. 


(iii) However, since $X$ is not required to be continuous in $\o$, we are not able to prove the desired regularity of $Y$. In particular, we are not able to prove that $Y$ and $\tilde Y$ are $\dbP_0$-indistinguishable. Consequently, we cannot verify rigorously that $Y$ solves RBSDE \reff{RBSDE}.  

(iv) In Theorem \ref{thm-optimal}, although $\hat X$ is also not continuous in $\o$, due to its special structure we proved in \cite{ETZ0} that the Snell envelope $Y$ has certain regularity in $\o$, which is crucial for proving the optimality of $\t^*$ in Theorem \ref{thm-optimal}.
\qed}\end{rem}
 
 We now establish the partial comparison principle.
\begin{prop}
\label{prop-semi-partial}
Let Assumption \ref{assum-F} hold and $\si=I_d$.   Let $ u^1\in \Usub$ be a viscosity subsolution of PPDE (\ref{PPDE}) and $ u^2\in\overline C^{1,2}(\L)$ satisfying $\cL u^2(t,\o) \ge 0$ for all $(t,\o) \in [0, T)\times \O$.  If $ u^1(T,\cd) \le u^2(T,\cd)$, then  $ u^1\le u^2$ on $\L$. 

The result also holds if we assume instead  that $u^1\in\overline C^{1,2}(\L)$ satisfies $\cL u^1(t,\o) \le 0$ for all $(t,\o) \in [0, T)\times \O$ and $u^2\in\Usup$ is a viscosity supersolution of PPDE (\ref{PPDE}).
\end{prop}
\proof As in Proposition \ref{prop-comparison}, we shall only prove $u^1_0 \le u^2_0$ under the additional condition \reff{Gmonotone}. Let $u^2 \in \overline C^{1,2}(\L)$ with corresponding  $\ch_i$, $i\ge 0$. Assume to the contrary that
$
c :=\frac{1}{2T}[ u^1_0 - u^2_0] > 0
$,
and denote 
\beaa
X_t := (u^1_t - u^2_t)^+ + ct.
\eeaa
By Definition \ref{defn-barC} (iii), $X$ satisfies the requirements in Remark \ref{rem-RBSDE}. Let $(\tilde Y, \tilde Z, \tilde K)$ and $\t^*$ be defined as in Remark \ref{rem-RBSDE} (i).   Then one can easily see that $\tilde K_t = 0$ for $t < \t^*$, and thus
\beaa
2cT = X_0 \le \tilde Y_0 = \overline \cE^L_0[\tilde Y_{\t^*}] =  \overline \cE^L_0[X_{\t^*}].
\eeaa 
This implies that
\beaa
 \overline \cE^L_0[X_T] = cT  < 2cT \le  \overline \cE^L_0[X_{\t^*}],
 \eeaa
and thus  $\dbP_0[\t^* < T] > 0$.  On the other hand,  apply Remark \ref{rem-RBSDE} (ii) to $\t^*$, then there exists $\o^*\in \O$ such that $t^* := \t^*(\o^*) < T$ and $Y_{t^*}(\o^*)= \tilde Y_{t^*}(\o^*)$. Thus
\beaa
X_{t^*}(\o^*) =\tilde Y_{t^*}(\o^*) = Y_{t^*}(\o^*) = \sup_{\t \in \cT^{t^*}} \overline\cE^L_{t^*}[ X^{t^*,\o^*}_\t].
\eeaa
In particular, this implies that 
\beaa
(u^1-u^2)^+(t^*, \o^*) + ct^* =  X_{t^*}(\o^*) \ge  \overline\cE^L_{t^*}[ X^{t^*,\o^*}_T] = cT,
 \eeaa
and thus $(u^1-u^2)(t^*,\o^*) > 0$. Assume without loss of generality that $t^* \in [\ch_i(\o^*), \ch_{i+1}(\o^*))$, and we may choose the $\ch\in \cH^{t^*}$ in Definition \ref{defn-barC} (ii)  small enough so that  $(u^1-u^2)^{t^*,\o^*}> 0$ on $[t^*, \ch)$.   Now  following the arguments in Proposition \ref{prop-comparison}, in particular by replacing the $(u^2)^{t^*,\o^*}$ there with the $\tilde u^{t^*,\o^*}\in C^{1,2}(\L^{t^*})$ in Definition \ref{defn-barC} (ii),   we can easily obtain the desired contradiction.
\qed

We now turn to comparison and uniqueness. First, define
 \bea\label{baru}
 \overline u(t,\o) 
 := 
 \inf\big\{\psi(t,{\bf 0}): 
           \psi\in\overline\cD(t,\o)\big\},
 \q 
 \underline u(t,\o) 
 := 
 \sup\big\{\psi(t,{\bf 0}): \psi\in  \underline\cD(t,\o)\big\},
\eea
where, for the $\cL$ in \reff{semiPPDE2} and denoting by $\cL^{t,\o}$ the corresponding operator on the shifted space with coefficient $F^{t,\o}$,
 \bea
 \label{cD}
 \left.\ba{lll}
 \dis \overline\cD(t,\o) 
 :=
 \Big\{ \psi \in \overline C^{1,2}(\L^t):
        ~\cL^{t,\o} \psi \ge 0~\mbox{on}~[t,T)\times\O^t
        ~\mbox{and}~
        \psi_T \ge \xi^{t,\o}\Big\},
 \\
 \dis \underline\cD(t,\o) 
 :=
 \Big\{\psi \in \overline C^{1,2}(\L^t):
       ~\cL^{t,\o}\psi \le 0~\mbox{on}~[t,T)\times\O^t
       ~\mbox{and}~
       \psi_T \le \xi^{t,\o}\Big\}.
\ea\right.
\eea
Following the arguments in the consistency Theorem  \ref{thm-consistency}, one can easily show that
\bea
\label{u<u}
\underline u &\le& \overline u.
\eea
A crucial step for our proof is to show that equality holds in the last inequality.

\begin{prop}\label{prop-perron}
Let Assumption \ref{assum-F} hold with $\si = I_d$, and $F$ is uniformly continuous in $(t,\o)$ under $\dbf_\infty$. Then have $\overline u = \underline u$. 
\end{prop}

We then have the following wellposedness result. 

\begin{thm}\label{thm-wellposedness}
Assume all the conditions in Proposition \ref{prop-perron} hold true.
\\
{\rm (i)}  Let $u^1\in \Usub$ be a viscosity subsolution and $u^2\in \Usup$ a viscosity supersolution of semilinear PPDE (\ref{semiPPDE2}), in the sense of Definition \ref{defn-semi-viscosity}, with $u^1_T \le \xi\le  u^2_T$. Then  $u^1\le u^2$ on $\L$.
\\
{\rm (ii)} The semilinear PPDE \reff{semiPPDE2} with terminal condition $\xi$ has a unique viscosity solution $u\in UC_b(\L)$, in the sense of Definition \ref{defn-semi-viscosity}.
\end{thm}
\proof First by the partial comparison principle Proposition  \ref{prop-semi-partial},  we have $u^1 \le \overline u$ and $\underline u \le u^2$. Then  Proposition \ref{prop-perron} implies $u^1\le u^2$ immediately, which implies further the uniqueness of viscosity solution. Finally by Proposition \ref{prop-semi} we have the existence. 
\qed

\vspace{5mm}

\no{\bf Proof of Proposition \ref{prop-perron}.} Without loss of generality, we shall only prove $\overline u(0, {\bf 0}) \le \underline u(0, {\bf 0})$.  In light of Proposition \ref{prop-change}, we may also assume without loss of generality that 
\bea
\label{semi-fmonotone}
F
~\mbox{is nonincreasing in $y$.}
\eea
For any $\e>0$, we denote 
 \beaa
 \left.\ba{lll}
 O_\e:=\{x\in\dbR^d: |x|<\e\},~~ \overline O_\e :=\{x\in\dbR^d: |x|\le\e\},~~ \pa O_\e :=\{x\in\dbR^d: |x|=\e\};\\
   Q^{\e}_t :=  [t,T)\times O_\e,\q \overline  Q^{\e}_t 
  := [t, T]\times \overline O_\e,\q
 \pa  Q^{\e}_t 
 :=\big([t,T]\times\partial O_\e\big) 
   \cup \big(\{T\}\times O_\e\big).
   \ea\right.
 \eeaa
 Moreover, we introduce the following space of discrete sequences:
 \beaa
 \Pi^\e_n &:=& \Big\{\pi_n = (t_i, x_i)_{0\le i\le n}: t_0 = 0, x_0 = {\bf 0}, t_i < t_{i+1}\wedge T ~\mbox{and}~|x_i|\le \e,~\mbox{for all}~i\Big\}.
 \eeaa
Now for any $\pi_n \in \Pi^\e_n$, and any $(t,x) \in  Q^\e_{t_n}$, define
 \beaa
 \ch^{t,x, \e}_1
 &:=& 
 T\wedge\inf\{s\ge t: |B^t_s + x| =\e\},\\
 \ch^{t,x, \e}_{i+1}
 &:=& 
 T\wedge\inf\big\{s \ge \ch^{t,x,\e}_{i}: 
                  |B^t_s - B^t_{\ch^{t,x,\e}_{i}}| = \e\big\}, 
 ~i\ge 1.
 \eeaa
We denote by $\hat B^{\e,\pi_n, t, x}(\o)$ the linear interpolation  $(t_i, \sum_{j=0}^i x_j)_{0\le i\le n}$ and $(\ch^{t,x, \e}_i(\o), \sum_{j=0}^n x_j + x + B^t_{\ch^{t,x,\e}_i}(\o))_{i\ge 1}$, namely, abbreviating $\ch_i := \ch_i^{t,x,\e}$, 
\beaa
\hat B^{\e,\pi_n, t, x}_s&=& \sum_{j=0}^i x_j +\frac{s-t_i}{t_{i+1}-t_i}x_{i+1},\q 0\leq i\leq n-1, s\in [t_i,t_{i+1}];\\
\hat B^{\e,\pi_n, t, x}_s&=&\sum_{j=0}^n x_j + x + B^t_{\ch_i} +\frac{s-t_i}{t_{i+1}-t_i}(B^t_{\ch_{i+1}}-B^t_{\ch_{i}}),\q i\ge 1, s\in [\ch_i,\ch_{i+1}].
\eeaa

Define 
\beaa
\th^\e_n\big(\pi_n; (t,x)\big) := \cY^{\e,\pi_n, t, x}_t
\eeaa
where, denoting $\ch^{t,x,\e}_{0} := t$ and omitting the superscripts $^{\e,\pi_n, t, x}$,
\beaa
\cY_s 
= 
\xi(\hat B) 
+ \int_s^T F\Big(r,\sum_{i\ge 0}
                   \hat B_{\cd\wedge \ch^{t,x,\e}_{i}}
                   \1_{[\ch^{t,x,\e}_{i}, \ch^{t,x,\e}_{i+1})}(r),
                 \cY_r, \cZ_r
            \Big)dr - \int_s^T \cZ_r\cdot dB_r,~\dbP_0^t\mbox{-a.s.}
\eeaa
We remark that 
\beaa
F\Big(r,\sum_{i\ge 0}
                   \hat B_{\cd\wedge \ch^{t,x,\e}_{i}}
                   \1_{[\ch^{t,x,\e}_{i}, \ch^{t,x,\e}_{i+1})}(r),
                 \cY_r, \cZ_r
            \Big) = \sum_{i\ge 0} F\Big(r,
                   \hat B_{\cd\wedge \ch^{t,x,\e}_{i}},
                 \cY_r, \cZ_r
            \Big) \1_{[\ch^{t,x,\e}_{i}, \ch^{t,x,\e}_{i+1})}(r)
            \eeaa
            is well defined and $\dbF$-adapted. One can easily prove that the deterministic function $\th^\e_n:= \th^\e_n(\pi_n; \cd)  \in C^{1,2}( Q^\e_{t_n}) $ and that  $\th^\e_n$ is continuous on the boundary $\big((t_n, T)\times \pa O_\e \big) \cup \big(\{T\}\times \overline O_\e\big)$. Indeed, it satisfies the following standard PDE in $ Q^\e_{t_n}$:
\bea
\label{then}
&&-\pa_t \th^\e_n - \frac12 I_d: D^2 \th^\e_n
- F(s, \o^{\pi_n, (T, {\bf 0})}, \th^\e_n, D \th^\e_n ) = 0 ~\mbox{in}~ Q^\e_{t_n}, 
\eea
with boundary conditions
\beaa
\th^\e_n(T,x) = \xi(\o^{\pi_n, (T,x)}),~|x|\le \e;\q  \th^\e_n(t,x) = \th^\e_{n+1}(\pi_n, (t,x); t, {\bf 0}),~t\in (t_n, T), x\in \pa O_\e.
\eeaa
where $\o^{\pi_n, (T,x)}$ denotes the linear interpolation of $(t_i, \sum_{j=0}^i x_j)_{0\le i\le n}$, $(T, \sum_{j=0}^n x_j + x)$, and is deterministic.

We now let $\ch^\e_i := \ch^{0,{\bf 0},\e}_i$, and $\hat B^{\e}$ the linear interpolation of  $\{ (\ch^{\e}_i, B_{\ch^{\e}_i})_{i\ge 0}\}$. Define
\beaa
\psi^\e(t,\o) := \sum_{n=0}^\infty \th^\e_n\big((0, {\bf 0}), (\ch^{\e}_i, B_{\ch^{\e}_i}- B_{\ch^\e_{i-1}} )_{1\le i\le n}; t, B_t - B_{\ch^\e_n}\big) \1_{[\ch^\e_n, \ch^\e_{n+1})}(t).
\eeaa
One can easily check that $\ch^\e_i$ satisfies Definition \ref{defn-barC} (i) and $\psi^\e$ satisfies  Definition \ref{defn-barC} (iii). Moreover,  for each $n\ge 0$, $\o\in \O$, and $\ch_n^\e(\o) \le t < \ch^\e_{n+1}(\o)$, we have $(t,\o_t- \o_{\ch_n^\e(\o)}) \in  Q^\e_{\ch_n^\e(\o)}$. Set $(t_0, x_0) := (0, {\bf 0})$, $(t_j, x_j) := (\ch_{j}^\e(\o), \o_{\ch_{j}^\e(\o)}-\o_{\ch_{j-1}^\e(\o)})$, $j=1,\cds, n$, and  let $\d >0$ be small enough such that $(t,\o_t- \o_{t_n}) \in  Q^{\e,\d}_{t_n} :=  [t_n, T-\d) \times O_{\e-\d}$.  One may modify $\th^\e_n$ outside of $  Q^{\e,\d}_{t_n}$ to obtain $\tilde \th^\e_n((t_j, x_j)_{0\le j\le n}; \cd) \in C^{1,2}([t_n, T]\times \dbR^d)$. Now by setting $\ch := \inf\{s \ge t: (s,  B^t_s + \o_t - \o_{t_n}) \notin  Q^{\e,\d}_{t_n} \} < (\ch_{i+1}^\e)^{t,\o}$ and $\tilde{ (\psi^\e)}^{t,\o}(s, B^t)  := \tilde \th^\e_n((t_j, x_j)_{0\le j\le n}; s, B^t_s +  \o_t - \o_{t_n})$,  it is clear that $\psi^\e$ satisfies Definition \ref{defn-barC} (ii). That is, $\psi^\e\in \overline C^{1,2}(\L)$ with corresponding hitting times $\ch^\e_n$.

One may easily check further that  $\psi^\e(T,\o) = \xi(\hat B^\e)$, and
\bea
\label{Lpsie}
 -\pa_t \psi^\e - \frac12 I_d:\pa^2_{\o\o} \psi^\e 
 - F\Big(s, \sum_{i\ge 0}
            \hat B^\e_{\cd\wedge \ch^{\e}_{i}}
            \1_{[\ch^{\e}_{i}, \ch^{\e}_{i+1})}, 
         \psi^\e, \pa_\o\psi^\e
    \Big) = 0.
 \eea
Notice that $\|\hat B^\e - B\|_T \le 2\e$. Then
 \bea
 \label{fge}
 \|\xi(\hat B^\e)- \xi(B)\|
 \le 
 \rho_0(2\e),\q \Big|F\Big(s,  \sum_{i=0}^\infty\hat B^\e_{\cd\wedge \ch^{\e}_{i}} \1_{[\ch^{\e}_{i}, \ch^{\e}_{i+1})} , y,z\Big) - F(s,   B , y,z)\Big|\le \rho_0(2\e).
 \eea
Set
 \beaa
 \overline\psi^\e := \psi^\e + \rho_0(2\e)[1 + T- t],\q    \underline\psi^\e := \psi^\e - \rho_0(2\e)[1 + T-t].
 \eeaa
 Then 
 \beaa
&& \overline\psi^\e \ge \psi^\e,\q  \overline\psi^\e \in \overline C^{1,2}(\L),\q  \overline\psi^\e(T,\o) \ge \psi^\e(T,\o) + \rho_0(2\e) = \xi(\hat B^\e) + \rho_0(2\e)\ge \xi(B),
\eeaa
and,  by \reff{semi-fmonotone}, \reff{fge}, and \reff{Lpsie}
\beaa
&&   -\pa_t \overline\psi^\e - \frac12 I_d:\pa^2_{\o\o} \overline \psi^\e -  F(s,   B_\cd , \overline \psi^\e, \pa_\o\overline\psi^\e) \\
&\ge&-\pa_t \psi^\e +  \rho_0(2\e) - \frac12 I_d:\pa^2_{\o\o} \psi^\e  - F(s,B_\cd , \psi^\e, \pa_\o\psi^\e)\\
&\ge&
-\pa_t \psi^\e  - \frac12 I_d:\pa^2_{\o\o} \psi^\e
- F\Big(s,\sum_{i\ge 0}\hat B^\e_{\cd\wedge \ch^{\e}_{i}}
                        \1_{[\ch^{\e}_{i}, \ch^{\e}_{i+1})},  
        \psi^\e, \pa_\o\psi^\e
   \Big) =0.
\eeaa
That is, $\overline\psi^\e \in \overline\cD(0,{\bf 0})$. Then $\overline u(0, {\bf 0}) \le \overline\psi^\e(0, {\bf 0})$. Similarly, one can prove   $\underline u(0, {\bf 0}) \ge \underline\psi^\e(0, {\bf 0})$. Thus
\beaa
\overline u(0, {\bf 0}) - \underline u(0, {\bf 0}) \le \overline\psi^\e(0, {\bf 0}) -  \underline\psi^\e(0, {\bf 0}) = 2\rho_0(2\e)(1 + T).
\eeaa
Send $\e\to 0$, we obtain $\overline u(0, {\bf 0}) \le \underline u(0, {\bf 0})$. This, together with \reff{u<u}, implies that $\overline u(0, {\bf 0}) = \underline u(0, {\bf 0})$.
\qed

We shall remark that  the regularity of $\th^\e_n$ is quite subtle. In fact,   in general $\th^\e_n$ may be discontinuous on $\{t_n\}\times \pa O_\e$.  However, we do not need the continuity at those points in above proof. 

\begin{rem}
\label{rem-comparison}
{\rm The above proof of Proposition \ref{prop-perron} takes advantage of  the following three facts in the semi-linear case, which do not hold anymore in the fully nonlinear case in \cite{ETZ2}.

(i) The proof of partial comparison principle Proposition  \ref{prop-semi-partial} uses the RBSDE theory, which applies implicitly the dominated convergence theorem. In the fully nonlinear case, in order to avoid the application of the dominated convergence theorem, we shall modify the space $\overline C^{1,2}(\L)$ slightly.

(ii) The functions $\th^\e_n$ can be defined via BSDEs. In the fully nonlinear case, in particular when there is no representation formula, we shall prove the existence of $\th^\e_n$ satisfying \reff{then}  in an abstract way in \cite{ETZ2}.

(iii) The functions $\th^\e_n$ are already in $C^{1,2}(Q^\e_{t_n})$. In the fully nonlinear case, this is typically not true, and then we shall approximate $\th^\e_n$ by smooth functions.
\qed}
\end{rem}

\section{First order PPDEs}
\label{sect-First}
\setcounter{equation}{0}

In this section we study the following first order PPDE:
\bea\label{equation-fst}
\{-\pa_t u -G(.,u,\pa_\o u)\}(t,\o)=0, &&(t,\o)\in \L, 
\eea
where $G: \L\times \dbR \times \dbR^d \rightarrow \dbR$ verifies the following counterpart of Assumption \ref{assum-G}:

\begin{assum}\label{assum-G-fst}

{\rm (i)}  $G(\cd, y,z)\in \dbL^0(\L)$ for any fixed $(y,z)$, and $|G(\cd, 0, {\bf 0})|\le C_0$.
\\
{\rm (ii)} $G$ is  uniformly Lipschitz continuous in $(y,z)$ with a Lipschitz constant $L_0$,  and locally uniformly continuous in $(t,\o)$ under $\dbf_\infty$, namely for any $(t,\o)$, there exists $\rho_{t,\o}$ such that
\beaa
\sup_{y, z}|G(\tilde t, \tilde \o, y, z)-G(t, \o, y, z)| \le \rho_{t,\o}(\dbf_\infty((t,\o), (\tilde t, \tilde \o))).
\eeaa
\end{assum}
We note that here we require $G$ to be locally uniformly continuous in $(t,\o)$, which is stronger than Assumption \ref{assum-G} (iv), but weaker than the uniform continuity required  in Theorem \ref{thm-wellposedness} or in \cite{ETZ2}. The uniform regularity  is used mainly for the proof of comparison principle. However, in this case we will employ some compactness arguments, and then continuity implies locally uniform continuity.

This PPDE was studied by Lukoyanov \cite{Lukoyanov} by using the compactness of the set $\Omega^t_L$ defined below. In this section we explain briefly how to reduce  our general formulation to this case so as to adapt to Lukoyanov's arguments. However, we emphasize again  that this type of compactness argument encounters a fundamental difficulty in  the second order case, see Remark \ref{rem-Holder} below. We shall establish the wellposedness of second order equations in our accompanying paper \cite{ETZ2} by using the optimal stopping result Theorem \ref{thm-optimal}.

As is well known, a first order HJB equation corresponds to deterministic control. Similar to Section \ref{sect-semilinear}, in this case we may restrict our probability measures to degenerate ones with $\beta = {\bf 0}$, see Remark \ref{rem-semilinear}. We then define for any $L>0$:
\bea
\label{first-cP}
\cP_L^t := \{ \dbP^\a: \a: [t, T] \to \dbR^d, |\a|\le L\} &\mbox{where}& dB^t_s = \a_s ds,~ \dbP^\a\mbox{-a.s.}
\eea
and the corresponding nonlinear expectations $\overline\cE^L_t$, $\underline\cE^L_t$, and nonlinear optimal stopping problems $\overline\cS^L_t$, $\underline\cS^L_t$, etc.  in an obvious way.  Denote
\bea
\label{OL}
\O^t_L := \{\o \in \O^t: \o ~\mbox{is Lipschitz continuous with Lipschitz constant $L$}\},
\eea
and $ \L^t_L := [t, T]\times \O^t_L$, $\cP^t_\infty := \cup_{L>0} \cP^t_L$, $\O^t_\infty := \cup_{L>0} \O^t_L$, $\L^t_\infty := \cup_{L>0} \L^t_L$.
As in \cite{Lukoyanov}, one can easily check that
\bea
\label{OLproperty}
\left.\ba{c}
\dbP(\O^t_L) = 1~\mbox{for all}~\dbP\in \cP^t_L,\q \O^t_L~\mbox{is compact and}~\O^t_\infty \subset\O~\mbox{is dense  under $\|\cd\|_T$};\\
\mbox{for}~s<t, \o\in \O^s_L, \tilde\o\in \O^t_L,~\mbox{we have}~\o\otimes_t\tilde\o \in \O^s_L.
\ea\right.
\eea

\begin{rem}
\label{rem-Holder}
{\rm All the above properties are important in Lukoyanov's approach for first order PPDEs, especially for proving the comparison principle. In the second order case, for example for the semilinear PPDEs considered in Section \ref{sect-semilinear}, since $\dbP_0^t(\O^t_\infty) = 0$, the set $\O^t_\infty$ is not appropriate.  One may consider to enlarge the space: for $0<\a<1$ and $L>0$, let
\beaa
\O^t_{\a,L} := \{\o \in \O^t: \o ~\mbox{is H\"{o}lder-$\a$ continuous  with H\"{o}lder constant $L$}\},~~\O^t_{\a,\infty} := \cup_{L>0} \O^t_{\a,L}.
\eeaa
Then for $\a<\frac12$ we have
\beaa
\dbP_0^t(\O^t_{\a,\infty}) = 1,\q \O^t_{\a, L}~\mbox{is compact and}~\O^t_{\a,\infty} \subset\O~\mbox{is dense  under $\|\cd\|_T$}
\eeaa
However, the last property in \reff{OLproperty} fails in this case:
\beaa
\mbox{for}~s<t, \o\in \O^s_{\a,L}, \tilde\o\in \O^t_{\a,L},~\mbox{in general }~\o\otimes_t\tilde\o \notin \O^s_{\a,L}.
\eeaa
This is the main reason why we were unable to extend this approach to second order case.
\qed}
\end{rem}
Notice that PPDE \reff{equation-fst} involves only derivatives $\pa_t u$ and $\pa_\o u$, we thus introduce the following definitions.

\begin{defn}\label{defn-spaceC11} 
We say a process $u\in C^0(\L^t)$ is  in $C^{1,1}(\L^t )$ if there exist $\pa_t u \in C^0(\L^t)$  and $\pa_\o u \in C^0(\L^t, \dbR^d)$ such that,
\bea
\label{Ito_fst}
d u_s = \pa_t u_s  ds+ \pa_\o u_s \cd d B^t_s ,
~~t\le s\le T, ~\dbP\mbox{-a.s. for all}~\dbP\in \cP^t_\infty.
\eea
\end{defn}
It is obvious that $\pa_t u$ and $\pa_\o u$, if they exist, are unique on $\L^t_\infty$. Then, since $\O^t_\infty\subset \O$ is dense and $\pa_t u$, $\pa_\o u$ are continuous, we see that they are unique in $\L^t$. 

For all $u\in \dbL^0(\L)$, $(t,\o) \in \L$ with $t<T$, and  $L>0$, define
 \bea\label{cA_fst}
 \ba{c}
 \underline\cA^{\!L}\!u(t,\o) 
 := 
 \Big\{\f\in C^{1,1}(\L^{\!t}):
       (\f-u^{t,\o})_t 
       = 0=
      \underline\cS^L_t\big[(\f-u^{t,\o})_{\cd\wedge \ch}
                       \big] ~\mbox{for some}~\ch\in \cH^t
\Big\},
 \\
 \overline\cA^{\!L}\!u(t,\o) 
:= 
\Big\{\f \in C^{1,1}(\L^{\!t}):
      (\f-u^{t,\o})_t
      =0=
      \overline\cS^L_t\big[(\f-u^{t,\o})_{\cd\wedge \ch}
                       \big]  ~\mbox{for some}~\ch\in \cH^t
 \Big\}.
 \ea
 \eea
We then define viscosity solutions exactly as in Definition \ref{defn-viscosity}.  We may easily check that all the results in this paper, when reduced to first order PPDEs, still hold under this new definition. In particular, the examples in Section \ref{sect-eg-first} are still valid, and the value function of the deterministic control problem is a viscosity solution to the corresponding first order path dependent  HJB equation. 

We remark that our Definition \ref{defn-spaceC11}  of derivatives is equivalent to Lukoyanov's notion of derivatives, which is defined via Taylor expansion. Moreover, instead of using nonlinear expectation as  in \reff{cA_fst},  Lukoyanov uses test functions $\f$ such that $\f - u$ attains pathwise local maximum (or minimum) at $(t,\o)$. So, modulus some minor technical difference, in spirit a viscosity solution (resp. subsolution, supersolution) in our sense is equivalent to a viscosity solution (resp. subsolution, supersolution) in Lukoyanov's sense. Indeed, our following comparison principle and uniqueness result for first order PPDEs follows from almost the same arguments as that of \cite{Lukoyanov}. We nevertheless sketch a proof for completeness. 

\begin{thm}
\label{thm-first}
Let Assumption \ref{assum-G-fst} hold true. Let $u^1\in \Usub$ (resp. $u^2\in \Usup$) be a viscosity subsolution (resp. supersolution) of PPDE \reff{equation-fst}, in the sense of Definition \ref{defn-viscosity} and modified in the context of this section. If $u^1(T,\cd) \le u^2(T,\cd)$ on $\O$, then $u^1 \le u^2$ on $\L$.
\end{thm}
\proof Let $u^1$ (resp. $u^2$) be a viscosity $L$-subsolution (resp. $L$-supersolution)  for some $L\ge L_0$. Assume by contradiction that $c_0 := u^1_0 - u^2_0 >0$.
For $\e>0$, define 
\beaa
&\Phi_\e(t, \o; \tilde t, \tilde \o) := u^1(t,\o) - u^2(\tilde t, \tilde \o) - {c_0\over 4T}[2T - t - \tilde t] - {1\over \e} \Psi_\e(t, \o; \tilde t, \tilde \o)&\\
&\mbox{where} \q\Psi_\e(t, \o; \tilde t, \tilde \o):= |t-\tilde t|^2+ |\o_t - \tilde\o_{\tilde t}|^2 + \int_0^T |\o_{t\wedge r} - \tilde \o_{\tilde t\wedge r}|^2dr.&\nonumber
\eeaa
Then
$
c_\e:=\sup_{(t, \o; \tilde t, \tilde \o) \in (\L_L)^2} \Phi_\e(t, \o; \tilde t, \tilde \o)  \ge  \Phi_\e(0, {\bf 0}; 0, {\bf 0})  = {c_0\over 2} >0.
$
By compactness of the space, there exists $(t_\e, \o^\e; \tilde t_\e, \tilde \o^\e)\in (\L_L)^2$ such that $\Phi_\e (t_\e, \o^\e; \tilde t_\e, \tilde \o^\e) = c_\e$.  
Note that $u^1$ is bounded from above and $u^2$ is bounded from below, then one can easily see that
\beaa
\Psi_\e (t_\e, \o^\e; \tilde t_\e, \tilde \o^\e)  \le C\e,&\mbox{which implies}& \lim_{\e\to 0}\dbf_\infty(t_\e, \o^\e; \tilde t_\e, \tilde \o^\e)=0.
\eeaa
Moreover, assuming without loss of generality that $t_\e \le \tilde t_\e$ and since $\Phi_\e (t_\e, \o^\e; \tilde t_\e, \tilde \o^\e) \ge \Phi_\e (\tilde t_\e, \tilde \o^\e; \tilde t_\e, \tilde \o^\e)$, it follows from \reff{USC} that 
\beaa
0\le {1\over \e} \Psi_\e (t_\e, \o^\e; \tilde t_\e, \tilde \o^\e) \le u^1(t_\e, \o^\e) - u^1(\tilde t_\e, \tilde \o^\e) - {c_0\over 4T}[\tilde t_\e-t_\e] \to 0.
\eeaa
Now, if $\tilde t_\e = T$, then $u^2(\tilde t_\e, \tilde \o^\e) \ge u^1(\tilde t_\e, \tilde \o^\e)$ and thus
\beaa
{c_0\over 2} \le  \Phi_\e (\tilde t_\e, \tilde \o^\e; \tilde t_\e, \tilde \o^\e) \le u^1(t_\e,\o^\e) - u^1(\tilde t_\e, \tilde \o) - {c_0\over 4T}[\tilde t_\e - t_\e] - {1\over \e} \Psi_\e(t, \o; \tilde t, \tilde \o)\to 0.
\eeaa
This is a contradiction. Thus we have $t_\e \le \tilde t_\e < T$ (or $\tilde t_\e \le t_\e <T$) when $\e$ is small enough.

We now define test functions:
\beaa
\f_1(t,\o) &:=& u^2(\tilde t_\e, \tilde \o^\e) + {c_0\over 4T}[2T - t - \tilde t_\e] + {1\over \e} \Psi_\e(t, \o; \tilde t_\e, \tilde \o^\e) - c_\e,\\
 \f_2(\tilde t,\tilde \o) &:=& u^1( t_\e, \o^\e) - {c_0\over 4T}[2T - t_\e - \tilde t] - {1\over \e} \Psi_\e(t_\e, \o^\e; \tilde t, \tilde \o) + c_\e.
\eeaa
It is straightforward to check that $\f_1 \in \underline \cA_L u^1(t_\e, \o^\e)$, $\f_2 \in \overline \cA_L u^2(\tilde t_\e, \tilde \o^\e)$, and 
\beaa
\pa_t \f_1(t_\e, \o^\e) =-{c_0\over 4T} + {2\over \e}[t_\e-\tilde t_\e],\q \pa_\o \f_1(t_\e, \o^\e) ={2\over \e}\Big[ [\o^\e_{t_\e} - \tilde \o^\e_{\tilde t_\e}] + \int_{t_\e}^T [\o^\e_{t_\e} - \tilde \o^\e_{\tilde t_\e\wedge r}]dr\Big];\\
 \pa_t \f_2(\tilde t_\e, \tilde \o^\e) ={c_0\over 4T} - {2\over \e}[\tilde t_\e- t_\e],\q  \pa_\o \f_2(\tilde t_\e, \tilde \o^\e)  =-{2\over \e}\Big[ [ \tilde \o^\e_{\tilde t_\e}-\o^\e_{t_\e}] + \int_{\tilde t_\e}^T [ \tilde \o^\e_{\tilde t_\e}] -  \o^\e_{ t_\e\wedge r}]dr\Big].
\eeaa
Note that $0 < \Phi_\e (\tilde t_\e, \tilde \o^\e; \tilde t_\e, \tilde \o^\e) \le u^1(t_\e,\o^\e) - u^2(\tilde t_\e, \tilde \o)$.  As standard, we may assume without loss of generality that $G$ is decreasing in $y$. Then it follows from the viscosity property of $u^1$, $u^2$ that
\beaa
{c_0\over 2T} &=&  \pa_t \f_2(\tilde t_\e, \tilde \o^\e) - \pa_t \f_1(t_\e, \o^\e) \le G(\cd, u_2, \pa_\o \f_2)(\tilde t_\e, \tilde \o^\e) - G(\cd, u_1, \pa_\o \f_1)(t_\e, \o^\e) \\
&\le& G\Big(\tilde t_\e, \tilde \o^\e, u_1(t_\e, \o^\e), \pa_\o \f_2(\tilde t_\e, \tilde \o^\e)\Big) - G\Big(t_\e, \o^\e, u_1(t_\e, \o^\e), \pa_\o \f_1(t_\e, \o^\e)\Big)\\
&\le& \sup_{y, z} |G(\tilde t_\e, \tilde \o^\e, y,z) - G(t_\e, \o^\e,y, z)| + L_0 |\pa_\o \f_2(\tilde t_\e, \tilde \o^\e)- \pa_\o \f_1(t_\e, \o^\e)|.
\eeaa
By the compactness of $\O_L$, $\{(t_\e, \o^\e), \e>0\}$ has a limit point $(t^*, \o^*)\in \L_L$, and we may assume without loss of generality that  $\dis\lim_{\e\to 0}\dbf_\infty ((t_\e, \o^\e), (t^*, \o^*))=0$. Then it follows from the locally uniform continuity of $G$ that $\dis\lim_{\e\to 0} \sup_{y, z} |G(\tilde t_\e, \tilde \o^\e, y,z) - G(t_\e, \o^\e,y, z)|=0$.
Moreover, 
\beaa
&& |\pa_\o \f_2(\tilde t_\e, \tilde \o^\e)- \pa_\o \f_1(t_\e, \o^\e)| = {2\over \e}\Big|\int_{t_\e \wedge \tilde t_\e}^{t_\e\vee \tilde t_\e} [\o^\e_{\t_\e\wedge r} - \tilde \o^\e_{\tilde t_\e \wedge r}]dr\Big|\\
 &\le&  {2\over \e}\Big[|t_\e - \tilde t_\e||\o^\e_{\t_\e} - \tilde \o^\e_{\tilde t_\e}| + \int_{t_\e \wedge \tilde t_\e}^{t_\e\vee \tilde t_\e}[|\o^\e_{t_\e\wedge r}-\o^\e_{t_\e}| + | \tilde \o^\e_{\tilde t_\e \wedge r}- \tilde \o^\e_{\tilde t_\e}|]dr\Big]\\
 &\le& {C\over \e}\Big[|t_\e - \tilde t_\e|^2 + |\o^\e_{\t_\e} - \tilde \o^\e_{\tilde t_\e}|^2 \Big]\to 0.
\eeaa
This implies ${c_0\over 2T} \le 0$, contradiction.
\qed



\end{document}